\documentclass{article}

\usepackage[utf8]{inputenc}
\usepackage{xspace}

\usepackage{amsmath,amsthm,amssymb}
\usepackage{mathtools}
\usepackage{altvars}

\usepackage[linesnumbered,ruled,vlined]{algorithm2e}
\SetKwInput{KwInput}{Input}                
\SetKwInput{KwOutput}{Output}              

\usepackage[hidelinks]{hyperref}
\usepackage[capitalise]{cleveref}
\crefname{equation}{}{}
\Crefname{equation}{Equation}{Equations}
\crefname{figure}{Fig.}{Figs.}

\usepackage{xcolor}
\definecolor{heteroppca}{RGB}{ 53, 96,240}
\definecolor{homoppca}{RGB}{255,160,253}
\definecolor{group1}  {RGB}{ 73,227,176}
\definecolor{group2}  {RGB}{255,166, 23}
\definecolor{inv}  {RGB}{ 10,112, 56}
\definecolor{sqinv}{RGB}{247,  5,  5}
\definecolor{heteropca}{RGB}{115,222,255}
\definecolor{heteroppcaknown}{RGB}{255,194, 38}

\usepackage[caption=false]{subfig}
\usepackage{graphicx}
\usepackage{tikz,pgfplots}
\pgfplotsset{compat=1.16}

\newcommand{\rootfind} {global maximization\xspace}
\newcommand{\heppca} {HePPCAT\xspace}

\newcommand{\diag}{\operatorname{diag}}
\newcommand{\tr}{\operatorname{tr}}

\newcommand{\xmath}[1] {\ensuremath{#1}\xspace}
\newcommand{\blmath}[1] {\xmath{\mathbf{#1}}}                        
\newcommand{\defequ} {\triangleq}
\newcommand{\paren}[1] {\xmath{\left( #1 \right)}}
\newcommand{\normfrob}[1] {\xmath{\left\| #1 \right\|_{\mathrm{F}}}} 
\newcommand{\normfrobs}[1] {\xmath{\| #1 \|_{\mathrm{F}}}}           

\newcommand{\alphatl} {\xmath{\alpha_{t,\ell}}}
\newcommand{\betatl} {\xmath{\beta_{t,\ell}}}
\newcommand{\zetatl} {\xmath{\zeta_{t,\ell}}}
\newcommand{\Btl} {\xmath{B_{t,\ell}}}
\newcommand{\ctl} {\xmath{c_{t,\ell}}}
\newcommand{\brLl} {\xmath{\breve{\mathcal{L}}_{\ell}}}
\newcommand{\betj} {\xmath{\beta_j}}
\newcommand{\gamj} {\xmath{\gamma_j}}
\newcommand{\clj} {\xmath{c_{\ell,j}}}
\let\nlalg\nl 
\renewcommand{\nl} {\xmath{n_{\ell}}}
\newcommand{\vtl} {\xmath{v_{t,\ell}}}
\newcommand{\vl} {\xmath{v_{\ell}}}
\newcommand{\sumjk} {\sum_{j=0}^k}
\newcommand{\sumjz}[1][] {\sum#1_{j \in \clJ_0}}
\newcommand{\sumjnoz}[1][] {\sum#1_{j \notin \clJ_0}}
\newcommand{\nolim}{\nolimits}
\newcommand{\nonum}{\nonumber}
\newcommand{\alf} {\xmath{\alpha}}

\DeclareMathOperator*{\argmax}{argmax}

\newcommand{\legendentry}[5]{
  \draw[#4,fill,opacity=0.25] (#1,#2-0.12) -- (#1+0.5,#2-0.12)
    -- (#1+0.5,#2+0.12) -- (#1,#2+0.12);
  \draw[#3,#4,very thick] (#1,#2) -- (#1+0.5,#2);
  \node[anchor=north west] at (#1+0.5,#2+0.225) {#5};
}

\makeatletter
\newenvironment{customlegend}[1][]{%
    \begingroup
    \let\addlegendimage=\pgfplots@addlegendimage
    \let\addlegendentry=\pgfplots@addlegendentry
    \pgfplots@init@cleared@structures
    \pgfplotsset{#1}%
}{%
    \pgfplots@createlegend
    \endgroup
}%
\makeatother

\graphicspath{{fig/}}

\usepackage{fullpage}
\newcommand{\oneortwocol}[2]{#1}
\newcommand{\IEEEPARstart}[2]{#1#2}
\newcommand{\blfootnote}[1]{
  \begingroup
  \renewcommand\thefootnote{}\footnote{#1}%
  \addtocounter{footnote}{-1}%
  \endgroup
}

\begin{document}

\title{\heppca: Probabilistic PCA for Data with Heteroscedastic Noise}

\author{%
  David Hong%
  \thanks{%
    Department of Statistics and Data Science,
    University of Pennsylvania,
    Philadelphia, PA, 19104 USA
    (email: dahong67@wharton.upenn.edu).
  }%
  \and
  Kyle Gilman%
  \thanks{%
    Department of Electrical Engineering and Computer Science,
    University of Michigan,
    Ann Arbor, MI, 48109 USA.
  }%
  \and
  Laura Balzano%
  \footnotemark[2]
  \and
  Jeffrey A. Fessler%
  \footnotemark[2]
}

\maketitle

\blfootnote{%
  D. Hong and L. Balzano were supported in part by
    ARO YIP award W911NF1910027.
  D. Hong was also supported in part by
    NSF BIGDATA grant IIS 1837992,
    the Dean's Fund for Postdoctoral Research of the Wharton School,
    and NSF Mathematical Sciences Postdoctoral Research Fellowship DMS 2103353.
  L. Balzano was also supported in part by
    NSF CCF-1845076 and the IAS Charles Simonyi Endowment.
  D. Hong, K. Gilman, L. Balzano and J. A. Fessler were supported in part by
    NSF BIGDATA IIS-1838179.
}
\blfootnote{DOI: \href{https://doi.org/10.1109/TSP.2021.3104979}{10.1109/TSP.2021.3104979}}
\blfootnote{%
  \textcopyright\, 2021 IEEE.
  Personal use of this material is permitted.
  Permission from IEEE must be obtained for all other uses, in any current or future media,
  including reprinting/republishing this material for advertising or promotional purposes,
  creating new collective works, for resale or redistribution to servers or lists,
  or reuse of any copyrighted component of this work in other works.
}

\begin{abstract}
%
Principal component analysis (PCA) is a classical and ubiquitous method
for reducing data dimensionality,
but it is suboptimal for heterogeneous data
that are increasingly common in modern applications.
PCA treats all samples uniformly
so degrades when the noise is heteroscedastic across samples,
as occurs, e.g., when samples come from sources of heterogeneous quality.
This paper develops a probabilistic PCA variant
that estimates and accounts for this heterogeneity
by incorporating it in the statistical model.
Unlike in the homoscedastic setting,
the resulting nonconvex optimization problem
is not seemingly solved by singular value decomposition.
This paper develops a
heteroscedastic probabilistic PCA technique (\heppca)
that
uses efficient alternating maximization algorithms to
jointly estimate both the underlying factors
and the unknown noise variances.
Simulation experiments illustrate
the comparative speed of the algorithms,
the benefit of accounting for heteroscedasticity,
and the seemingly favorable optimization landscape
of this problem.
Real data experiments on environmental air quality data show that
\heppca
can give a better PCA estimate than techniques that do not account for heteroscedasticity.

\end{abstract}


\section{Introduction}
\label{sec:intro}

\IEEEPARstart{P}{rincipal} component analysis (PCA)
is a workhorse method for unsupervised dimensionality reduction.
It plays a foundational role in the analysis of modern high-dimensional data,
and continues to be successfully applied across all of engineering and science.
However, PCA does not account for samples having heterogeneous quality
and instead treats them uniformly.
Consequently, the performance of PCA can degrade dramatically
under heteroscedastic noise;
its ability to  discover underlying components
is sometimes essentially determined
by the noisiest samples alone \cite{hong2018apo}.

At the same time,
heterogeneous quality among samples is common in practice,
arising easily when samples are obtained under varying conditions.
For example, in the field of air quality monitoring,
there is a wide array of sensors available for different entities to deploy:
governments use very high-quality sensors
that require regular maintenance but are very accurate,
and individuals purchase off-the-shelf sensor devices
that can be deployed and left alone but have much less reliable output.
These devices are measuring the same phenomenon through very different noise characteristics.
In the field of analytical chemistry,
\cite{cochran1977swp} considers
spectrophotometric data
that are averages over increasingly long windows of time.
This heterogeneity
arises naturally when measuring signals
that undergo both periods of rapid change (requiring short windows)
as well as periods of relatively stable behavior (allowing for longer windows).
The shorter windows cannot denoise by averaging as much,
resulting in heteroscedasticity.
Another source of heteroscedasticity
is changing ambient conditions;
e.g., \cite{tamuz2005cse} considers astronomical data
with atmospheric noise that varies across nights.
As large datasets are increasingly formed
by combining samples from diverse sources,
one can expect that heteroscedastic noise
will be the norm.
Modern data analysis needs PCA methods
that are robust to heterogeneity
and make effective use of all the available data.

This paper develops a heteroscedastic probabilistic PCA technique (\heppca)
that attains robustness to heteroscedastic noise
by incorporating it in the statistical likelihood.
The method
jointly estimates
both the latent factors
as well as the unknown sample-wise noise variances.
Additionally, if a block of samples are expected
to have equal noise variance
(e.g., because they are from the same source or sensor),
the proposed approach seamlessly incorporates this knowledge
and can yield significantly improved estimates.
A further extension to the case where some variances are known
and some are unknown is straightforward.
Unlike the homoscedastic setting,
the resulting optimization problem
seems not to have a direct SVD solution.
Because it is nonconvex and nontrivial,
we develop and compare several alternating ascent algorithms.

\heppca is an extension of our previous work~\cite{hong2019ppf}
that considered data with \emph{known} heterogenous noise variances
and focused on estimating the latent factors alone.
In this paper, the noise variances are \emph{unknown and jointly estimated}
with the latent factors.
This extension is important in practice
because heterogeneous data often have unknown noise variances.
It is also nontrivial to do efficiently.
As discussed in \cref{sec:em},
the Expectation Maximization (EM)
approach used for the latent factors in \cite{hong2019ppf}
does not readily yield an efficient approach in this joint estimation setting.
Thus, we develop and study efficient block coordinate ascent algorithms
that alternate between updating estimates of the latent factors and
estimates of the noise variances.

\Cref{sec:ppca} describes the model and the resulting optimization problem for \heppca.
\Cref{sec:rel} discusses related works.
\Cref{sec:em} derives a natural EM approach,
and explains why the resulting M-step is challenging.
This difficulty motivates alternating approaches
that are derived in \cref{sec:alg}
and compared in \cref{sec:complexity,sec:exp:comp}.
\Cref{sec:exp:stat} carries out several experiments
illustrating the favorable statistical performance of \heppca.
\Cref{sec:exp:real_data} illustrates \heppca on real data.
\Cref{sec:exp:init} investigates
the seemingly favorable landscape of the nonconvex objective,
illustrating that the proposed algorithms
appear to converge from even random initializations.
A Julia package implementing \heppca and
code to reproduce all experiments
will be
available online at:
\url{https://gitlab.com/heppcat-group/heteroscedastic-probabilistic-pca}.


\section{Heteroscedastic Probabilistic PCA}
\label{sec:ppca}

As in
\cite{hong2019ppf},
we model $n_1 + \cdots + n_L = n$ data samples in $\bbR^d$
from $L$ noise level groups as:
\begin{equation}
  \bmy_{\ell,i} = \bmF \bmz_{\ell,i} + \bmvarepsilon_{\ell,i}
  ,
  \quad
  i \in \{1,\dots,n_\ell\},\
  \ell \in \{1,\dots,L\}
  ,
  \label{eq:ppca:generative_model}
\end{equation}
where
$\bmF \in \bbR^{d \times k}$ is a deterministic factor matrix to estimate,
$\bmz_{\ell,i} \sim \clN(\bm0_k,\bmI_k)$
  are independent and identically distributed (i.i.d.) coefficients,
$\bmvarepsilon_{\ell,i} \sim \clN(\bm0_d,v_\ell \bmI_d)$ are i.i.d. noise vectors,
and $v_1,\dots,v_L$ are deterministic noise variances to estimate.
Equivalently, the samples are independent with distributions
\begin{equation*}
  \bmy_{\ell,i} \sim \clN(\bm0_d, \bmF\bmF' + v_\ell \bmI_d)
  ,
\end{equation*}
and joint log-likelihood, dropping the $\ln(2\pi)^{-nd/2}$ constant:
\begin{align}  \label{eq:likelihood}
  \clL(\bmF,\bmv)
  &\defequ
  \frac{1}{2}
  \sum_{\ell=1}^L \Big[
    n_\ell \ln\det(\bmF\bmF' + v_\ell \bmI_d)^{-1}
  \oneortwocol{}{\\&\qquad\qquad\qquad\nonumber}
    - \tr\big\{\bmY_\ell'(\bmF\bmF' + v_\ell \bmI_d)^{-1}\bmY_\ell\big\}
  \Big]
  ,
\end{align}
where
$\bmY_\ell
\defequ [\bmy_{\ell,1},\dots,\bmy_{\ell,n_\ell}]
\in \bbR^{d \times n_\ell}$
for $\ell \in \{1,\dots,L\}$
are the sample matrices associated with each of the $L$ groups.
Note that $\bmF'$ denotes the matrix transpose for real-valued $\bmF$;
the methods generalize easily to complex matrices
using the Hermitian transpose.

\begin{figure}
    \centering
    \oneortwocol{
      \hfill
      \subfloat[$n_1 = 20, n_2 = 10$]{\includegraphics[trim=0 0 6.2cm 0,clip,scale=0.6]{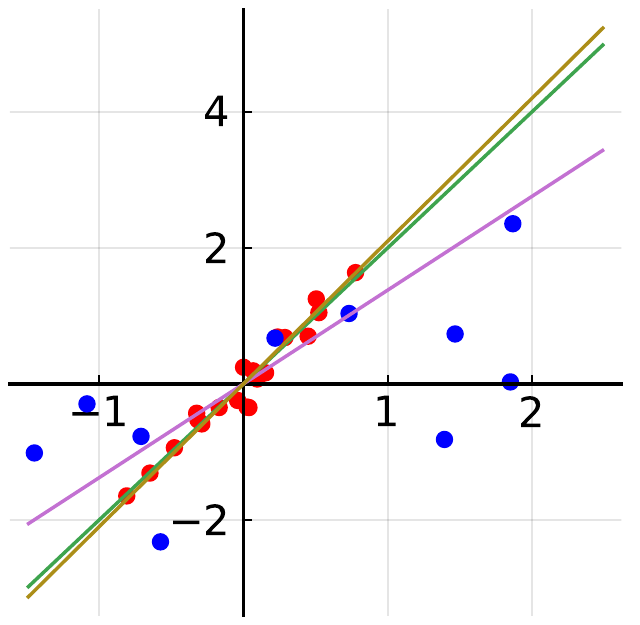}} 
      \hfill
      \subfloat[$n_1 = 10, n_2 = 20$]{\includegraphics[scale=0.6]{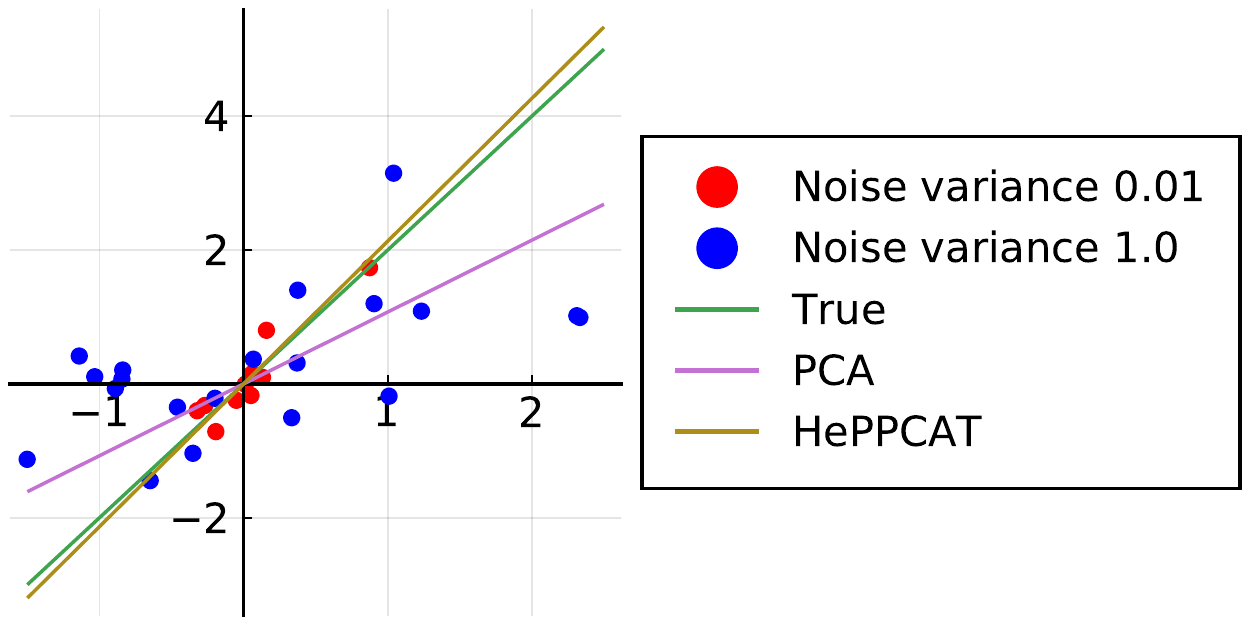}}
      \hfill
    }{
      \subfloat[$n_1 = 20, n_2 = 10$]{\includegraphics[trim = 0 0 6.2cm 0, clip, scale=0.4]{rank1hppca_1.pdf}} 
      \hspace{2mm}
      \subfloat[$n_1 = 10, n_2 = 20$]{\includegraphics[scale = 0.4]{rank1hppca_2.pdf}}
    }
    \caption{Illustrative heteroscedastic examples with $k=1$ factor
      and $L=2$ noise variances $v_1 = 0.01$ and $v_2 = 1$.
      \heppca estimates variances to account for heteroscedasticity,
      recovering the true latent subspace much better than PCA.}
    \label{fig:rank-1_visual}
\end{figure}

Given the sample matrices $\bmY_1,\dots,\bmY_L$
and the rank $k$,
\heppca estimates the latent factors $\bmF \in \bbR^{d \times k}$
and the noise variances $\bmv \defequ (v_1,\dots,v_L)$
by maximizing the statistical log-likelihood \cref{eq:likelihood}.
\Cref{fig:rank-1_visual}
shows an illustrative example
with $L=2$ noise variances $v_1 = 0.01$ and $v_2 = 1$.
When the noise is assumed homoscedastic, i.e., $L=1$,
this nonconvex optimization problem can be solved via eigendecomposition
of the sample covariance matrix \cite[Section 3.2]{tipping1999ppc},
but the same is not true in general.

The groupings give a natural way to
incorporate structural assumptions
by grouping together samples that are expected to have equal noise variance,
e.g., samples from the same source or sensor.
They are given and not estimated.
In the absence of such knowledge,
each sample can be given its own group
by taking $n_1 = \cdots = n_L = 1$ and $L = n$.
\heppca estimates a separate noise variance for each sample in that case,
and we study some of the resulting trade-offs in \cref{sec:stat:perf:blocks}.

Representing the factors by the rank-$k$ eigendecomposition
$\bmF\bmF' = \bmU \diag(\bmlambda) \bmU'$
where $\bmU = [\bmu_1,\dots,\bmu_k] \in \bbR^{d \times k}$
and $\bmlambda \in \bbR^k$
yields an alternative form for the likelihood:
\begin{align} \label{eq:likelihood:parts}
  &
  \clL(\bmU, \bmlambda, \bmv)
  \oneortwocol{}{\\&\qquad\nonumber}
  =
  \frac{1}{2}
  \sum_{\ell=1}^L \bigg[
    -
    n_\ell \bigg\{
      \sum_{j=1}^k \ln(\lambda_j + v_\ell) + (d-k) \ln v_\ell
    \bigg\}
  \oneortwocol{}{\\&\qquad\qquad\qquad\quad\nonumber}
    - \frac{\|\bmY_\ell\|_F^2}{v_\ell}
    +
    \tr\big\{
      \bmY_\ell'\bmU
      \bmW(\bmlambda,v_\ell)
      \bmU'\bmY_\ell
    \big\}
  \bigg]
  ,
\end{align}
with weighting matrices
\begin{equation} \label{eq:weighting:matrix}
  \bmW(\bmlambda,v)
  \defequ
  \diag\bigg(
    \frac{\lambda_1/v}{\lambda_1+v},
    \dots,
    \frac{\lambda_k/v}{\lambda_k+v}
  \bigg)
  .
\end{equation}
Maximizing \cref{eq:likelihood:parts} with respect to $\bmU$ resembles a weighted PCA,
but, unlike weighted PCA,
it is not readily solved by eigendecomposition in general
since the weight matrices $\bmW(\bmlambda,v_\ell)$ can vary with $\ell$.
Jointly optimizing further complicates the problem.
Following a review of related work,
the remainder of this paper investigates various alternating algorithms
for this joint maximization.

\section{Related works}
\label{sec:rel}

\subsection{Factor analysis and (homoscedastic) probabilistic PCA}
\label{sec:rel:factor:analysis}

In conventional factor analysis,
samples in $d$ dimensions are modeled as follows:%
\footnote{While more general versions exist, for simplicity, we omit the mean and focus here
on the conventional setting with Gaussian coefficients
and additive Gaussian noise
that is most closely related.}
\begin{equation*}
  \bmy_i = \bmF \bmz_i + \bmvarepsilon_i
  ,
  \quad
  i \in \{1,\dots,n\}
  ,
\end{equation*}
where $\bmF \in \bbR^{d \times k}$ contains the $k$ factors,
$\bmz_i \overset{iid}{\sim} \clN(\bm0_k,\bmI_k)$ are random coefficients,
and $\bmvarepsilon_i \overset{iid}{\sim} \clN(\bm0_d,\bmPsi)$ are random noise
with diagonal covariance $\bmPsi$.
The quantities
$\bmF, \bmz_i, \bmvarepsilon_i$ and $\bmPsi$ are all unknown.
Marginalizing out $\bmz_i$ and $\bmvarepsilon_i$
yields a model where only $\bmF$ and $\bmPsi$ are to be estimated
by maximizing the marginal likelihood.
In general, the column space of a maximum likelihood estimate for $\bmF$
will not coincide with the corresponding principal subspace of the data.
Indeed, factor analysis and PCA
are somewhat distinct approaches to dimensionality reduction;
see, e.g., \cite[Chapter 7]{jolliffe2002pca}.

However, maximum likelihood estimation does produce
the principal subspace
if the noise covariance is assumed to be isotropic,
i.e., $\bmPsi = v\, \bmI_d$ for some noise variance $v$.
This model is the setting of probabilistic PCA \cite{tipping1999ppc}.
In this case, the log-likelihood is
\begin{align*}
  \clL^\mathrm{PPCA}(\bmF,v)
  &\defequ
  \frac{1}{2}
  \Big[
    n \ln\det(\bmF\bmF' + v \bmI)^{-1}
  \oneortwocol{}{\\&\qquad\qquad\qquad\nonumber}
    - \tr\big\{\bmY'(\bmF\bmF' + v \bmI)^{-1}\bmY\big\}
  \Big]
  ,
\end{align*}
and is maximized by
$\bmF = \bhtU \diag^{1/2}(\htlambda_1-\brlambda, \dots, \htlambda_k-\brlambda)$
and $v = \brlambda$,
where the columns of $\bhtU \in \bbR^{d \times k}$
are principal eigenvectors of the sample covariance matrix
$(\bmy_1 \bmy_1' + \cdots + \bmy_n \bmy_n')/n$,
$\htlambda_1,\dots,\htlambda_k$ are the corresponding eigenvalues,
and $\brlambda$ is the average of the remaining $d-k$ eigenvalues
\cite{lawley1953modified,anderson1956sii,tipping1999ppc}.
Moreover, \cite{tipping1999ppc,roweis1998eaf} characterize stationary points
as well as the global maxima of the likelihood objective function.
They also derive an efficient expectation maximization (EM) algorithm
related to one derived for factor analysis \cite{rubin1982eaf},
and illustrate how the approach
naturally generalizes to similar models.

Here we develop a new probabilistic PCA method;
unlike previous settings,
the samples are no longer identically distributed.
Noise variances are now heterogeneous,
i.e., the noise is heteroscedastic across samples.
The resulting likelihood is no longer maximized
by scaled eigenvectors of the sample covariance,
so new algorithms are needed.
We developed
an EM algorithm for estimating the factors
given known noise variances
in \cite{hong2019ppf};
this paper extends that work
by jointly estimating both
the factors and the unknown noise variances.

\subsection{Accounting for heteroscedastic noise via weighted PCA}
\label{sec:rel:wpca}

A natural way to account for heteroscedastic noise
is to use a weighted PCA \cite[Section~14.2.1]{jolliffe2002pca}
that replaces the sample covariance with a weighted sample covariance.
Namely, given weights $w_1,\dots,w_L$,
weighted PCA returns the leading eigenvectors $\bhtu_1,\dots,\bhtu_k$
of $\sum_{\ell=1}^L w_\ell \bmY_\ell \bmY_\ell'$.
These eigenvectors solve the weighted optimization problem
\begin{equation*}
  \bhtU \defequ [\bhtu_1,\dots,\bhtu_k]
  \in
  \argmax_{\bmU \in \bbR^{d \times k} : \bmU'\bmU = \bmI_k}
  \sum_{\ell=1}^L w_\ell \normfrob{\bmU' \bmY_\ell}^2
  .
\end{equation*}
A typical choice for the weights is inverse noise variance,
i.e., $w_\ell = 1/v_\ell$, so
samples that are twice as noisy get half as much weight.
Doing so effectively whitens the noise,
is a type of maximum likelihood weighting \cite{young1941mle},
and can significantly improve performance \cite{hong2018owp:arxiv:v2}.
However, as analyzed in \cite{hong2018owp:arxiv:v2},
it can be better to more aggressively downweight noisier samples,
especially for low signal-to-noise ratio (SNR) regimes.
In particular,
optimal weights for recovery of any individual component
is between inverse noise variance
and square inverse noise variance,
which more aggressively downweights noisier samples.
Weighted PCA with general weights does not
have an obvious maximum likelihood formulation.

In contrast, this paper considers the maximum likelihood estimation
of underlying factors and noise variances jointly.
The resulting optimization problem does not appear to
reduce to PCA with a weighted sample covariance,
yielding a distinct approach to accounting for heteroscedasticity.

\subsection{Heteroscedasticity across features}

This paper focuses on noise that is heteroscedastic across samples,
i.e., the samples are of varying quality.
Noise can also be heteroscedastic \emph{across features}.
Indeed, much of the general literature on heteroscedasticity
focuses on this manifestation.
In the context of PCA,
recent works have begun to study
how to account for this form of heteroscedasticity.
Notably,
this heteroscedasticity
induces a bias along the diagonal of the covariance matrix,
causing conventional PCA to produce inaccurate components.
To correct for this bias,
\cite{zhang2018hpa:arxiv:v2} describes the HeteroPCA method
that treats the diagonal entries as missing and iteratively imputes the values.
Alternatively,
\cite{leeb2018oss:arxiv:v4,leeb2019mdf:arxiv:v3}
combine whitening of the data
with spectral shrinkages
tailored to optimize, e.g., matrix denoising.

\subsection{Accounting for heterogeneous clutter in RADAR}
\label{sec:rel:radar}

In the context of estimating low-rank clutter,
\cite{breloy2015cse,breloy2016rcm,sun2016lca}
model $n$ independent samples $\bmy_1,\dots,\bmy_n \in \bbC^d$
as
\begin{equation*}
  \bmy_i = \bmvarepsilon_i + \bmc_i
  , \quad
  i \in \{1,\dots,n\}
  ,
\end{equation*}
where $\bmvarepsilon_i \overset{iid}{\sim} \clC\clN(\bm0_d,\bmI_d)$ is complex white Gaussian noise,
and the clutter $\bmc_i \sim \clC\clN(\bm0_d,\tau_i \bmSigma)$ share a common rank-$k$ covariance $\bmSigma$
that is scaled by heterogeneous power factors~$\tau_i$.
Equivalently, $\bmy_i \sim \clC\clN(\bm0_d,\tau_i \bmSigma + \bmI_d)$ for $i \in \{1,\dots,n\}$.
The goal is to estimate $\tau_1,\dots,\tau_n$
and $\bmSigma$.

The low-rank covariance term $\tau_i \bmSigma$ is heterogeneous,
while the noise is homogeneous.
In contrast, the low-rank factor covariance in this paper
is common among all the samples,
and instead the \emph{noise is heterogeneous}.
The two models are related
through an unknown heterogeneous rescaling
because
\begin{equation*}
  \frac{1}{\sqrt{\tau_i}} \bmy_i \sim \clC\clN\{\bm0_d, \bmSigma + (1/\tau_i) \bmI_d\}
  , \quad
  i \in \{1,\dots,n\}
  ,
\end{equation*}
\oneortwocol{has}{corresponds to} a common low-rank factor covariance $\bmSigma$
and heteroscedastic noise with variances $1/\tau_1,\dots,1/\tau_n$.
As a result,
the two problems share some common challenges and approaches.
Notably, the EM factor update (\cref{sec:alg:F:em})
is essentially the same (up to rescaling)
as \cite[Section III-B]{sun2016lca}.

Nevertheless, the problems remain distinct due to the difference
in how the unknown heterogeneity manifests.
For example,
heterogeneous power factors are only identifiable up to scale,
since any change in scale can be absorbed by $\bmSigma$.
The heterogeneous noise model we study does not have this scale ambiguity.
Moreover,
the likelihood for heterogeneous noise
as a function of noise variances
has a similar form as
that for heterogeneous power factors,
but with significant differences.
Notably, decreasing a noise variance to zero
sends the log-likelihood to $-\infty$
with \emph{unbounded curvature}
in the common case where the data are not perfectly fit.
As a result,
approaches well-designed for updating power factor estimates,
such as the minorizer \cite[Proposition 1]{sun2016lca},
cannot always be directly applied to update the noise variance estimates.
Different algorithms are needed.

In the context of heterogeneous power factors,
\cite{besson2016bfa} derives bounds on estimation performance
and \cite{abdallah2020bss} places priors on the clutter subspace.
These are also interesting directions for future work
on heterogeneous noise.

\subsection{Matrix Factorization}
\label{sec:rel:opt}

The model \cref{eq:ppca:generative_model} that this paper focuses on
can also be interpreted
as a matrix factorization formulation
that has Gaussian coefficients and additive noise.
Within this framework there exist
generalizations where one assumes other coefficient and noise distributions
or even treats the factors or noise variances as random with a prior distribution.
That is, one may generalize \cref{eq:ppca:generative_model}
to allow other distributions on $\bmz_{\ell,i}, \bmvarepsilon_{\ell,i}$
and/or put a distribution on $\bmF$ or $v_1,\dots,v_L$.
There is a great deal of literature for factor analysis in a variety of settings,
such as non-negative matrix factorization \cite{gillis2017introduction},
Poisson matrix estimation \cite{cao2015poisson},
robust PCA \cite{lerman2018overview},
logistic PCA \cite{schein2003generalized},
and others \cite{bartholomew2011latent, udell2016generalized,shi2017survey, vaswani2018rethinking, liu2018pca}.
Extending these models for heterogeneous noise is an interesting direction for future work. 

In addition to this modeling work, great progress has been made in recent years to better understand why standard optimization algorithms perform well, even seeming to find the global minima/maxima,
when applied to nonconvex matrix factorization problems
\cite{chen2020noisy, ding2020leave, ge2017no, zhang2019sharp, ge2016matrix, bhojanapalli2016global, sun2016wan, park2017non}.
Three recent surveys summarize much of this progress
\cite{jain2017non, chi2019nonconvex, zhang2020symmetry}
and
thoroughly treat this related work.
An overview of minorize maximize (MM) techniques and how they are applied to related problems
can be found in \cite{sun2017mma}.
Recent guarantees applied specifically to EM are found in \cite{jain2017non}.
None of the existing results apply directly to our setting for two reasons:
first our model has noise that is not identically distributed across columns,
and second we seek to optimize over both the factor matrix $\bmF$
and the additive noise variances $\bmv$ in our maximum likelihood formulation.
If we consider only the problem of optimizing over the factor matrix,
one could potentially extend results from \cite{tipping1999ppc, ge2017no, bhojanapalli2016global}
to characterize the stationary points of our objective.

Numerous works that involve matrix factorization use a spectral initialization,
and several show that this initialization is sufficiently close to a good optima
\cite{chen2020noisy, lu2019phase}. 
We also use spectral initialization in our nonconvex optimization methods.


\section{Expectation Maximization}
\label{sec:em}

A natural way to maximize the log-likelihood \cref{eq:likelihood}
is through an expectation maximization (EM)
approach
that produces a sequence of iterates $\bmF_t$ and $\bmv_t$
with non-decreasing log-likelihood.
At each iteration,
an EM method
sets up a minorizer based on conditional expectation
(E-step)
that it then maximizes (M-step).
This section derives an EM minorizer for
the \heppca log-likelihood \cref{eq:likelihood}
at the iterates $\bmF_t$ and $\bmv_t$,
where $t$ denotes iteration index.
The resulting minorizer turns out to be challenging to maximize efficiently,
so instead \cref{sec:alg} proposes alternating algorithms,
where some of the updates are based on the EM minorizer derived here.

Taking as complete data
the samples $\bmY_1,\dots,\bmY_L$
and (unknown) coefficients $\bmZ_1,\dots,\bmZ_L$,
where
$\bmZ_\ell
\defequ [\bmz_{\ell,1},\dots,\bmz_{\ell,n_\ell}]
\in \bbR^{k \times n_\ell}$
for $\ell \in \{1,\dots,L\}$,
yields the following complete data log-likelihood
\begin{align}
  \oneortwocol{}{&}
  \clL_c(\bmF,\bmv)
  \oneortwocol{&}{}
  \defequ
  \ln p(\bmY, \bmZ; \bmF, \bmv)
  \oneortwocol{}{\nonumber \\ &\quad}
  = \ln p(\bmY | \bmZ; \bmF, \bmv) + \ln p(\bmZ; \bmF, \bmv)
  \nonumber \\
  &\oneortwocol{}{\quad}
  =
  \sum_{\ell=1}^L \bigg(
    - \frac{d n_\ell}{2} \ln v_\ell
    - \frac{\normfrob{\bmY_\ell - \bmF\bmZ_\ell}^2}{2v_\ell}
    - \frac{\normfrob{\bmZ_\ell}^2}{2}
  \bigg)
  , \label{eq:em:likelihood}
\end{align}
where \cref{eq:em:likelihood} drops
the constants $\ln(2\pi)^{-nd/2}$ and $\ln(2\pi)^{-nk/2}$.

For the E-step, take the expectation of \cref{eq:em:likelihood}
with respect to the conditionally independent distributions
(from Bayes' rule and the matrix inversion lemma):
\begin{equation} \label{eq:em:conddist}
  \bmz_{\ell,i} | \bmY, \bmF_t, \bmv_t
  \overset{\text{ind}}{\sim}
  \clN(
    \bmM_{t,\ell} \bmF_t' \bmy_{\ell,i},
    v_{t,\ell} \bmM_{t,\ell}
  )
  ,
\end{equation}
where $\bmM_{t,\ell} \defequ (\bmF_t'\bmF_t + v_{t,\ell}\bmI_k)^{-1}$,
yielding minorizer
\begin{align} \label{eq:em:estep}
  &\cbrL(\bmF, \bmv; \bmF_t, \bmv_t)
  \defequ
  \sum_{\ell=1}^L
    \bigg[
      - \frac{d n_\ell}{2} \ln v_\ell
      - \frac{\normfrob{\bmY_\ell}^2}{2 v_\ell}
  \oneortwocol{}{\\&\nonumber}
      + \frac{1}{v_\ell} \tr(\bmY_\ell' \bmF \bbrZ_{t,\ell})
      - \frac{1}{2 v_\ell}
      \tr\{
        \bmF'\bmF
        (\bbrZ_{t,\ell}\bbrZ_{t,\ell}' + n_\ell v_{t,\ell} \bmM_{t,\ell})
      \}
    \bigg]
  ,
\end{align}
where
$\bbrZ_{t,\ell}
\defequ \bmM_{t,\ell} \bmF_t' \bmY_\ell
\in \bbR^{k \times n_\ell}$
for $\ell \in \{1,\dots,L\}$,
and
\cref{eq:em:estep} drops terms
that are constant with respect to $\bmF$ and $\bmv$.

The corresponding M-step involves jointly maximizing \cref{eq:em:estep}
with respect to both $\bmv$ and $\bmF$,
but doing so is challenging
because the interaction of the variables remains complicated.
See \cref{sec:em:challenges} for more discussion.
However, optimization with respect to either (with the other fixed)
is relatively easy,
and \cref{sec:alg:F:em,sec:alg:v:em}
use this minorizer
to obtain efficient updates for the individual variables.


\section{Alternating Algorithms}
\label{sec:alg}

The challenge of jointly optimizing \bmF and \bmv
using \eqref{eq:likelihood}
or \eqref{eq:em:estep}
motivates approaches that alternate between:
a) optimizing \bmF for fixed \bmv, and
b) optimizing \bmv for fixed \bmF.
Namely, we consider a block-coordinate ascent of \cref{eq:likelihood}
with $\bmF$ and $\bmv$ as the two blocks of variables.
These sub-problems are simpler but
the sub-problem for updating \bmv
using \eqref{eq:likelihood}
or \eqref{eq:likelihood:parts}
is still nontrivial
so this section considers several methods
for updating $\bmv$ given $\bmF$.
When either the \bmF or \bmv update
involves the conditional expectation
with respect to some complete data,
then such alternation is an instance
of a space-alternating generalized EM (SAGE) algorithm
\cite{fessler1994sag}.


\subsection{Optimizing \texorpdfstring{$\bmF$}{F}
  for fixed \texorpdfstring{\bmv}{v} (via Expectation Maximization)}
\label{sec:alg:F:em}

Fixing $\bmv$ at $\bmv_t$,
maximizing the minorizer
$\cbrL(\bmF, \bmv_t; \bmF_t, \bmv_t)$
in \eqref{eq:em:estep}
with respect to $\bmF$
yields
the EM step of \cite{hong2019ppf}:
\begin{equation} \label{eq:alg:F:em}
  \bmF_{t+1}
  =
    \bigg(\sum_{\ell=1}^L \frac{\bmY_\ell \bbrZ_{t,\ell}'}{v_{t,\ell}} \bigg)
    \bigg(
      \sum_{\ell=1}^L
        \frac{\bbrZ_{t,\ell}\bbrZ_{t,\ell}'}{v_{t,\ell}}
        + n_\ell \bmM_{t,\ell}
    \bigg)^{-1}
  ,
\end{equation}
that we compute
via the SVD $\bmF_t = \bmU_t \bmLambda_t^{1/2} \bmV_t'$:
\begin{equation} \label{eq:alg:F:em:svd}
  \bmF_{t+1}
  =
    \bigg( \sum_{\ell=1}^L \frac{\bmY_\ell \btlZ_{t,\ell}'}{v_{t,\ell}} \bigg)
    \bigg(
      \sum_{\ell=1}^L \frac{\btlZ_{t,\ell}\btlZ_{t,\ell}'}{v_{t,\ell}}
      + n_\ell \bmD_{t,\ell}
    \bigg)^{-1}
    \bmV_t'
  ,
\end{equation}
where $\btlZ_{t,\ell} \defequ \bmD_{t,\ell} \bmLambda_t^{1/2} \bmU_t' \bmY_\ell$
and $\bmD_{t,\ell} \defequ (\bmLambda_t + v_{t,\ell} \bmI_k)^{-1}$
is easily inverted
because $\bmLambda_t$ and $\bmD_{t,\ell}$ are diagonal.
To show that this form is equivalent,
note that $\bmM_{t,\ell} = \bmV_t \bmD_{t,\ell} \bmV_t'$
and $\bbrZ_{t,\ell} = \bmV_t \btlZ_{t,\ell}$.
See \cite[Section 3]{hong2019ppf} and \cite[Section III-B]{sun2016lca} for similar derivations.


\subsection{Optimizing \texorpdfstring{\bmv}{v}
  for fixed \texorpdfstring{$\bmF$}{F}}
\label{sec:alg:v}

Fixing $\bmF$ at $\bmF_t$, maximization of \cref{eq:likelihood:parts}
with respect to $\bmv$ separates into $L$ univariate maximizations
(over $\vl \geq 0$)
of:
\begin{equation} \label{eq:nss:opt}
\clL_\ell(\vl)
\defequ
  -
  \sum_{j=0}^k \bigg\{
    \alpha_j \ln(\gamj + \vl)
    + \frac{\beta_j}{\gamj + \vl}
  \bigg\}
,
\end{equation}
where
$\alpha_0 \defequ d-k$,
$\beta_0 \defequ \normfrob{ (\bmI_d - \bmU_t\bmU_t')\bmY_\ell }^2/\nl$,
$\gamma_0 \defequ 0$,
\begin{align*}
  j &\geq 1
  : &
  \alpha_j &\defequ 1
  , &
  \beta_j &\defequ \|\bmY_\ell'\bmu_{t,j}\|_2^2/\nl
  , &
  \gamj &\defequ \lambda_{t,j}
  ,
\end{align*}
and $\bmU_t = [\bmu_{t,1},\dots,\bmu_{t,k}]$
and $\bmlambda_t = (\lambda_{t,1},\dots,\lambda_{t,k})$
are the eigenvectors and eigenvalues of $\bmF_t\bmF_t'$.
\Cref{eq:nss:opt} drops all terms from \cref{eq:likelihood:parts}
that are constant with respect to \vl
as well as a factor of $\nl/2$,
and we define
\begin{equation}
\label{eq:v:betatl}
  \clL_\ell(0)
  \defequ
  \clL_\ell(0^+)
  =
  \begin{cases}
    +\infty
    , &
    \text{if } \betatl \defequ \sumjz \betj = 0
    , \\
    -\infty
    , &
    \text{otherwise}
    ,
  \end{cases}
\end{equation}
where
\(
\clJ_0 \defequ \{ j \, : \, \gamj = 0 \}
.\)
Note also that $\clL_\ell(+\infty) = -\infty$
and $\forall_{\vl \in (0,\infty)} \; \clL_\ell(\vl) < \infty$.
Lacking an analytical solution for the critical points of 
\eqref{eq:nss:opt}
when $k > 1$,
we next describe several iterative methods
for maximizing
$\clL_\ell(\vl)$.

\subsubsection{Global maximization via root-finding}
\label{sec:alg:v:roots}

If 
$\betatl = 0$,
then \cref{eq:nss:opt} is maximized by $\vl = 0$.
Otherwise, $\clL_\ell(0^+) = -\infty$
and global maxima occur only at critical points.
Differentiating \cref{eq:nss:opt} with respect to \vl
yields
\begin{equation} \label{eq:nss:opt:diff}
\dot\clL_\ell(\vl)
\defequ
\sum_{j=0}^k
\bigg\{
  - \frac{\alpha_j}{\gamj + \vl}
  + \frac{\beta_j}{(\gamj + \vl)^2}
\bigg\}
.
\end{equation}
An upper bound for nonnegative roots of \cref{eq:nss:opt:diff}
can be obtained from general root bounds for polynomials,
e.g., \cite{marden1949gop,hong1998bfa}.
We exploit the structure here to find a specialized bound.
The $k+1$ summands in \eqref{eq:nss:opt:diff} are, respectively,
positive to the left and negative to the right of $\beta_j/\alpha_j - \gamj$.
As a result,
\begin{align*}
  \dot\clL_\ell(\vl) > 0 &\text{ for }
  \vl < \vl^{\min} \defequ \min_j (\beta_j/\alpha_j - \gamj)
  , \\
  \dot\clL_\ell(\vl) < 0 &\text{ for }
  \vl > \vl^{\max} \defequ \max_j (\beta_j/\alpha_j - \gamj)
  ,
\end{align*}
so all nonnegative critical points occur in $[\vl^{\min}, \vl^{\max}] \cap [0,\infty)$
and can be found, e.g., via interval root-finding%
\footnote{We used the Julia package
\href{https://github.com/JuliaIntervals/IntervalRootFinding.jl}{IntervalRootFinding.jl}.}
\cite[Ch.~8]{moore2009iti}.
Choosing the best among these critical points
yields global maximizers.

This update maximally ascends the likelihood,
and is perhaps the most natural choice.
However, finding all the roots can be computationally expensive.
Moreover, it is unclear whether fully maximizing the likelihood
in this step is desirable
since this update occurs within a broader alternating maximization.
The current estimate of $\bmF$ may be far from optimal,
so fully optimizing \bmv
might slow convergence.
These reasons motivate alternative methods
that we derive next.

\subsubsection{Expectation Maximization}
\label{sec:alg:v:em}

Although jointly updating $\bmF$ and $\bmv$ using
\cref{eq:em:estep} is challenging,
it is fairly easy to update $\bmv$ when $\bmF$ is fixed.
Replacing $\bmF$ in
\cref{eq:em:estep}
with the current estimate $\bmF_t$
and simplifying
leads to the following minorizer
of \eqref{eq:likelihood:parts}
with respect to~\vl:
\begin{equation} \label{eq:v:em:estep}
\cbrL(\bmF_t, \bmv; \bmF_t, \bmv_t)
=
\sum_{\ell=1}^L
  \frac{\nl}{2}
  \bigg(
    - d \ln \vl
    - \frac{\rho_{t,\ell}}{\vl}
  \bigg)
,
\end{equation}
where
\begin{align} \label{eq:v:em:residual}
\rho_{t,\ell}
&\defequ
\frac{1}{\nl}
\big[
  \normfrob{\bmY_\ell}^2
  - 2 \tr(\bmY_\ell' \bmF_t \bbrZ_{t,\ell})
\oneortwocol{}{\\ &\quad \hspace*{2em} \nonumber}
  +
  \tr \{ \bmF_t'\bmF_t
    (\bbrZ_{t,\ell}\bbrZ_{t,\ell}' + \nl \vtl \bmM_{t,\ell})
  \}
\big]
\\
&=
\normfrob{ (\bmI_d - \bmF_t \bmM_{t,\ell} \bmF_t') \bmY_\ell }^2 / \nl
+
\vtl
\tr( \bmF_t \bmM_{t,\ell} \bmF_t' )
. \nonumber
\end{align}
Maximizing \cref{eq:v:em:estep}
w.r.t.~\vl
leads to the
simple update:
\begin{equation}
v_{t+1,\ell} = \frac{\rho_{t,\ell}}{d}
.
\end{equation}
Since
$\bmF_t \bmM_{t,\ell} \bmF_t'
=\bmU_t \bmLambda_t (\bmLambda_t + \vtl \bmI_k)^{-1} \bmU_t'$,
expanding and simplifying yields the following alternative formula for \cref{eq:v:em:residual}:
\begin{equation} \label{eq:alg:v:em:svd}
  \rho_{t,\ell}
  =
  \sum_{j=0}^k \bigg(1 - \frac{\gamma_j}{\gamma_j+v_{t,\ell}}\bigg)^2 \beta_j
  +
  \vtl
  \sum_{j=1}^k \frac{\lambda_{t,j}}{\lambda_{t,j}+v_{t,\ell}}
  ,
\end{equation}
providing a more efficient form as well as a link to \cref{eq:nss:opt}.

\subsubsection{Difference of concave approach}
\label{sec:alg:v:diffconc}

The univariate objective
\cref{eq:nss:opt}
is a ``difference of concave'' or concave+convex cost function.
One standard way to optimize such functions
is to minorize each convex term with an affine function,
leading to the following concave minorizer
(ignoring constants):
\begin{equation}
\label{eq:v:doc}
  \ctlL_\ell(\vl; \vtl)
  \defequ
  -
  \sumjk \bigg\{
    \frac{\alpha_j}{\gamj + \vtl} \vl
    + \frac{\betj}{\gamj + \vl}
  \bigg\}
  .
\end{equation}
Concavity of \cref{eq:v:doc} eases maximization.
If its derivative,
\begin{equation*}
  \dot\ctlL_\ell(\vl; \vtl)
  \defequ
  \sumjk \bigg\{
    - \frac{\alpha_j}{\gamj + \vtl} + \frac{\betj}{(\gamj + \vl)^2}
  \bigg\}
  ,
\end{equation*}
is nonpositive at the origin,
i.e., $\dot\ctlL_\ell(0^+; \vtl) \leq 0$,
then \cref{eq:v:doc} is maximized by $\vl = 0$.
Otherwise, at least one $\beta_j > 0$
so \cref{eq:v:doc} is necessarily \emph{strictly} concave
and is maximized at its unique critical point over $\vl > 0$.
This critical point can be efficiently computed,
e.g., via bisection by noting that
$\dot\ctlL_\ell(\vl; \vtl) < 0$
for $\vl > \max_j\{\sqrt{(\betj/\alpha_j)(\gamj + \vtl)} - \gamj\}$.

\subsubsection{Quadratic solvable minorizer}
\label{sec:alg:v:quad}

To derive a MM approach with a simple update,
we separate the summation in
\eqref{eq:nss:opt}
into the terms where \gamj is zero and nonzero
and apply the affine minorizer
of \eqref{eq:v:doc}
to the $\ln$ terms where $\gamj > 0$
as follows:
\begin{align}
\brLl(\vl; \vtl) =&
-\alphatl \ln \vl - \frac{\betatl}{\vl}
- \zetatl \, \vl
\oneortwocol{}{\nonumber\\&}
- \sumjnoz \frac{\betj}{\gamj + \vl}
,
\label{eq:v-doc-split}
\end{align}
where
\(
\alphatl \defequ \sumjz \alf_j
,\)
\(
\zetatl \defequ \sumjnoz \frac{\alpha_j}{\gamj + \vtl}
,\)
and \betatl was defined in
\eqref{eq:v:betatl}.

For $j \notin \clJ_0$,
let
\(
\pi_j \defequ \frac{\gamj}{\gamj + \vtl}
\in (0,1)
\)
and rewrite
\eqref{eq:v-doc-split}
as
\begin{align}
& \brLl(\vl;\vtl)
=
-\alphatl \ln \vl - \frac{\betatl}{\vl} - \zetatl \, \vl
\oneortwocol{}{\nonum\\&\hspace*{3em}}
- \sumjnoz \frac{\beta_j}{\pi_j \paren{\frac{\gamj}{\pi_j}} + (1-\pi_j) \paren{\frac{\vl}{1 - \pi_j}}}
\nonum\\
\geq &
\ \phi(\vl;\vtl) \defequ
- \alphatl \ln \vl - \frac{\betatl}{\vl}
- \zetatl \, \vl
\oneortwocol{}{\nonum\\&\hspace*{2em}}
- \sumjnoz \beta_j \paren{\pi_j \frac{1}{\gamj/\pi_j} + (1-\pi_j) \frac{1}{\vl / (1-\pi_j)} }
\nonum\\&=
- \alphatl \ln \vl - \frac{\Btl}{\vl} - \zetatl \vl
,\end{align}
using the concavity of the function $-1/x$,
ignoring irrelevant constants,
and defining
\begin{align*}
\Btl &=
\betatl + \sumjnoz \betj \frac{\vtl^2}{(\gamj + \vtl)^2}
.\end{align*}
One can verify that, by design,
$\phi(\vtl;\vtl) = \brLl(\vtl; \vtl)$.
The choice of $\pi_j$ originates from an EM algorithm for PET
\cite{fessler1994sag,depierro:95:ame}.
Differentiating the concave minorizer
$\phi$ yields:
\begin{equation}
\label{eq:v,quad,solve}
0 = -\frac{\alphatl}{\vl} + \frac{\Btl}{\vl^2} - \zetatl
.\end{equation}
This equation is solvable by the quadratic formula
for
\(
\zetatl \vl^2 + \alphatl \vl - \Btl
\)
that has exactly one positive root.

\subsubsection{Cubic solvable minorizer}
\label{sec:alg:v:cubic}

Because $\gamj > 0$ in the final term
of the concave minorizer in \eqref{eq:v-doc-split},
that term has bounded curvature
for $\vl \geq 0$,
with maximum (absolute) curvature
\[
\clj = -2 \betj / \gamj^3
.\]
Thus we have
the following partially quadratic concave minorizer for
\eqref{eq:nss:opt}
(ignoring constants):
\begin{align}
& Q_l(\vl ; \vtl) =
-\alphatl \ln \vl - \frac{\betatl}{\vl} - \zetatl \, \vl
\oneortwocol{}{\nonumber\\&\hspace*{2em}}
+ \sumjnoz \bigg\{
\frac{\betj}{(\gamj + \vtl)^2} \vl
+ \frac{1}{2} \clj (\vl - \vtl)^2
\bigg\}
.\label{eq:v,quad,minor}
\end{align}
Differentiating and equating to zero yields
\begin{align}
0 &= \frac{-\alphatl}{\vl} + \frac{\betatl}{\vl^2}
+ \gamma_{t,\ell} + \ctl \, (\vl - \vtl)
,
\label{eq:v,cubic}
\\
\gamma_{t,\ell} &\defequ -\zetatl
+ \sumjnoz[\nolim]  \frac{\betj}{(\gamj + \vtl)^2}
, \nonumber\\
\ctl &\defequ \sumjnoz[\nolim] \clj
\nonumber
.
\end{align}
This \vl update corresponds
to finding the appropriate root
of a cubic polynomial.
One could apply multiple \vl updates
based on \eqref{eq:v,cubic}.
The fixed points of the resulting MM iteration
are identical to the roots described by
\eqref{eq:nss:opt:diff},
so this approach is essentially
an iterative root finding method
with the nice MM property of monotonically increasing
the log-likelihood.


\subsection{Convergence and stopping criterion}
\label{sec:conv}

All of the updates described above
for \bmF and \bmv
are based on minorizers,
so like all block MM methods
they provide updates
that ensure the log-likelihood is monotonically non-decreasing.
However, monotonicity alone is insufficient to ensure convergence
when the individual updates
may not have unique maximizers
\cite{powell:73:osd}.
An alternative to simple alternation
between the \bmF and \bmv blocks
is the ``maximum improvement'' variant
that calculates an update for both blocks
and chooses the one that increases the likelihood the most
\cite{chen:12:mbi}.
This variant ensures convergence
under modest regularity conditions
(not requiring convexity or uniqueness)
appropriate for the \heppca problem
\cite[Thm.~3]{razaviyayn:13:auc}.
To save computation,
we used the simpler alternating maximization approach
for the empirical results shown below.

A natural choice for stopping criterion
is to stop once the change in the factor matrix
is sufficiently small.
Namely, iterate until
$\normfrob{\bmF_{t+1} - \bmF_t} / \normfrob{\bmF_t} \leq \epsilon$,
where $\epsilon \geq 0$ is a user-provided tolerance,
as shown in \cref{alg:heppcat}.
That said, there are certainly other natural choices.
For example,
one could require sufficiently small changes
in the noise variance estimates
or in the log-likelihood.


\subsection{Initialization by homoscedastic PPCA}
\label{sec:init}

Without prior knowledge of the noise variances,
a natural choice to initialize $\bmv$ and $\bmF = \bmU \diag^{1/2}(\bmlambda)$
is the \emph{homoscedastic} PPCA solution
\cite[Section 3.2]{tipping1999ppc}:
\begin{align*}
  \bmU_0 &\defequ \bhtU
  , &
  \bmlambda_0 &\defequ (\htlambda_1-\brlambda, \dots, \htlambda_k-\brlambda)
  , &
  \bmv_0 &\defequ \brlambda \bm1_L
  ,
\end{align*}
where the $k$ columns of $\bhtU \in \bbR^{d \times k}$
are principal eigenvectors of the sample covariance matrix
$(\bmY_1 \bmY_1' + \cdots + \bmY_L \bmY_L')/n$,
$\htlambda_1,\dots,\htlambda_k$ are the corresponding eigenvalues,
and $\brlambda$ is the average of the remaining $d-k$ eigenvalues.

The \heppca optimization problem is nonconvex,
so better maximizers \emph{might} be found
by taking the best among many random initializations,
but we have not so far encountered such a case;
see, e.g., the experiments in \cref{sec:exp:init}.
The landscape of the objective
appears to be favorable despite its nonconvexity.
Moreover, initializing via homoscedastic PPCA
provides a reasonable and nicely interpretable choice.
If the samples are in fact close to homoscedastic,
this initialization is likely close to optimal already.
Even if not, it provides a homoscedastic baseline
to improve upon via the alternating updates of \cref{sec:alg:F:em,sec:alg:v}.
All the updates
are non-descending,
so all iterates are guaranteed to
have likelihood no worse than homoscedastic PPCA.

\begin{algorithm}
\def\nl{\nlalg} 
\DontPrintSemicolon
\SetAlgoHangIndent{0pt}
  
  \KwInput{Maximum number of iterations $T$, rank $k$, tolerance $\epsilon \geq 0$}
  \KwOutput{$\bmF \in \bbR^{d \times k}$, $\bmv \in \bbR_+^L$}
  \KwData{$[\bmY_1,\hdots,\bmY_L]$, $\bmY_\ell \in \bbR^{d \times n_\ell} \quad \ell=1,\hdots,L$}
  Initialize $\bmF_0$ and $\bmv_0$ via homoscedastic PPCA or random initialization \;
   \While{\textup{iterations} $t < T$}
   {
   	Update $\bmF_{t+1}$, fixing $\bmv$ at $\bmv_t$, using \eqref{eq:alg:F:em:svd}.\;
   	Update $\bmv_{t+1}$, fixing $\bmF$ at $\bmF_{t+1}$, using Root Finding
   	or one step of one of the following updates:
   	\begin{itemize}
    \item Expectation Maximization \eqref{eq:alg:v:em:svd}
    \item Difference of concave minorizer \eqref{eq:v:doc}
    \item Quadratic solvable minorizer \eqref{eq:v,quad,solve}
    \item Cubic solvable minorizer \eqref{eq:v,cubic}
    \end{itemize}\;
    
    \uIf{ $\normfrob{\bmF_{t+1} - \bmF_t} / \normfrob{\bmF_t} \leq \epsilon$ }{
    stop \;
    }
    $t \leftarrow t + 1$ \;
   }

\caption{HePPCAT}
\label{alg:heppcat}
\end{algorithm}


\section{Computational complexity}
\label{sec:complexity}

The primary sources of computational complexity
in the EM update \cref{eq:alg:F:em} for $\bmF$
are matrix multiplications and inverses.
For each $\ell \in \{1,\dots,L\}$,
computing $\bmM_{t,\ell}$ costs $O(k^2d + k^3)$
after which computing $\bbrZ_{t,\ell}$ costs $O(k^2d + kdn_\ell)$,
yielding a total cost of $O(L k^3 + L k^2d + kdn)$.
The remaining multiplications and inverses cost $O(kdn + k^2n + k^2d + k^3)$.
Combining these terms and noting that $k < d, n$ yields
$O(L k^2d + kdn)$.
The alternative form \cref{eq:alg:F:em:svd}
incurs an initial cost of $O(k^2d)$
to obtain the SVD of $\bmF_t$,
but gains efficiency
since $\bmD_{t,\ell}$ and $\btlZ_{t,\ell}$
then cost $O(k)$ and $O(kdn_\ell)$, respectively.
As a result, this form
has a final cost of $O(kdn)$ overall.

A leading order source of computational complexity
for all of the $v_\ell$ update methods is in calculating
the associated coefficients $\beta_0,\dots,\beta_k$.
Doing so incurs a cost of $O(kdn_\ell)$
for each $\ell \in \{1,\dots,L\}$,
yielding a cost of $O(kdn)$ overall.
For all the $v_\ell$ updates
(\cref{sec:alg:v:roots,sec:alg:v:em,sec:alg:v:diffconc,sec:alg:v:quad,sec:alg:v:cubic}),
the additional computational cost is independent of $d$ and $n$.
Thus, one might suppose (since $k \ll d,n$ typically)
that all the updates have essentially equal runtime.
However, this is \emph{not} the case.
Global maximization (\cref{sec:alg:v:roots})
and the difference of concave approach (\cref{sec:alg:v:diffconc})
both use iterative algorithms for root-finding,
and the runtime can depend significantly on not only $k$
but also properties of \cref{eq:nss:opt}.
Moreover, when $L$ is large,
e.g., for block sizes of $n_\ell = 1$,
constant factors not captured by computational complexity
can also have a significant impact.
\Cref{sec:exp:comp} compares the convergence speed
of the various updates in practice,
accounting for both their runtime costs
and per-iteration improvement in likelihood.

To further improve computational efficiency,
note that all the updates depend on $\bmY_\ell$
implicitly through $\bmY_\ell \bmY_\ell'$,
so one could replace $\bmY_\ell$ in the updates
with any proxy $\bbvY_\ell$
for which
$\bbvY_\ell \bbvY_\ell' = \bmY_\ell \bmY_\ell'$.
In some cases, e.g., when $n_\ell \gg d$,
doing so can yield significant savings.


\section{Comparison of update methods}
\label{sec:exp:comp}

\begin{figure*} \centering
  \subfloat[Convergence w.r.t. $\clL$. \label{fig:comp:hom:likelihood}]
    {\includegraphics[width=0.24\linewidth]{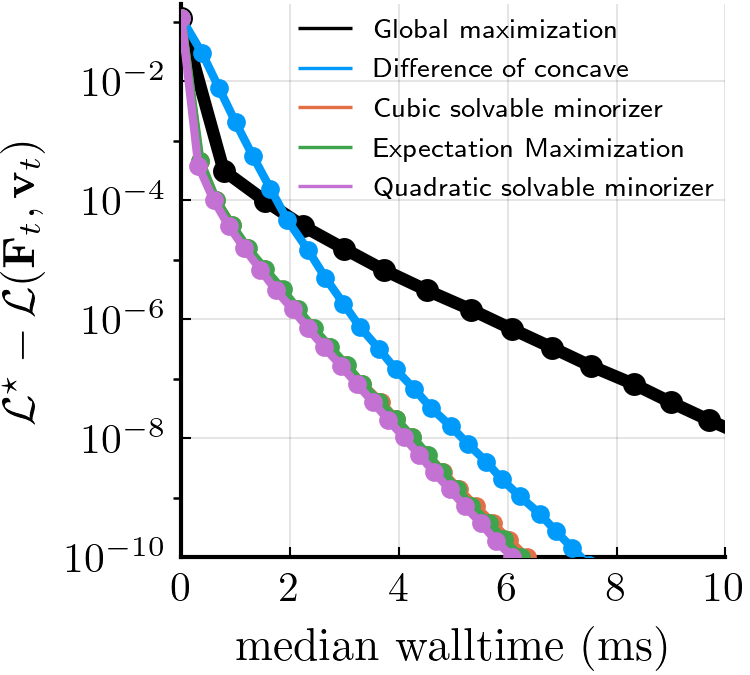}}%
  \hfill
  \subfloat[Convergence w.r.t. $\bmF$. \label{fig:comp:hom:factors}]
    {\includegraphics[width=0.24\linewidth]{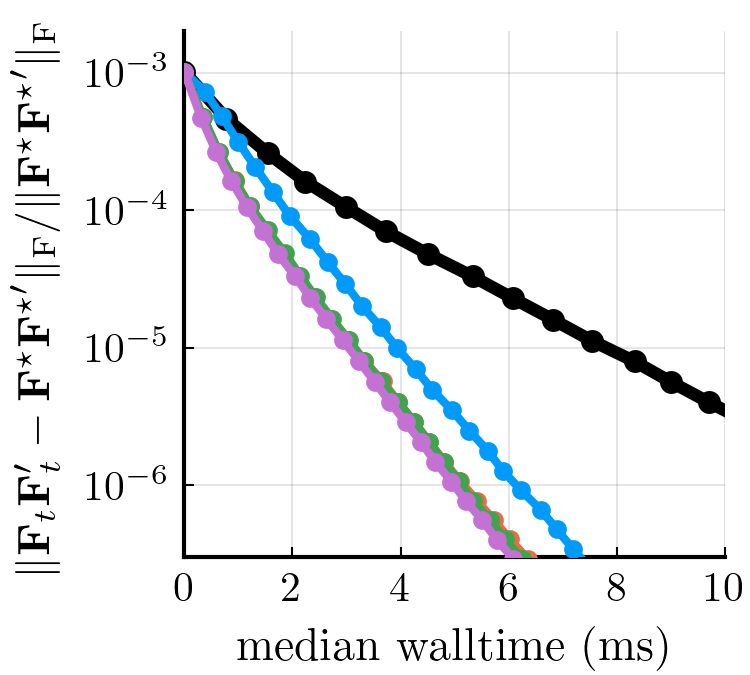}}%
  \hfill
  \subfloat[Minorizers w.r.t. $\bmv$ at homoscedastic PPCA initialization. \label{fig:comp:hom:minorizer}]
    {\includegraphics[width=0.24\linewidth]{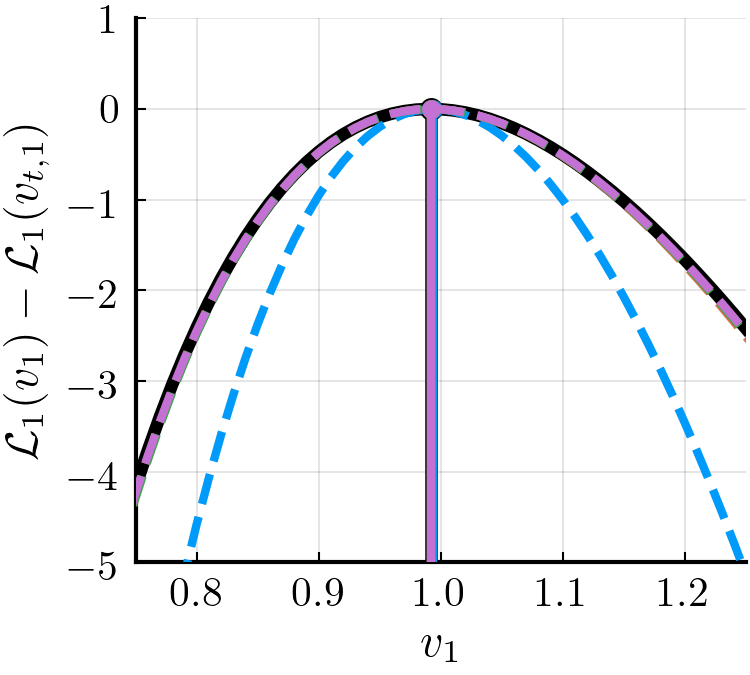}\;\;
     \includegraphics[width=0.24\linewidth]{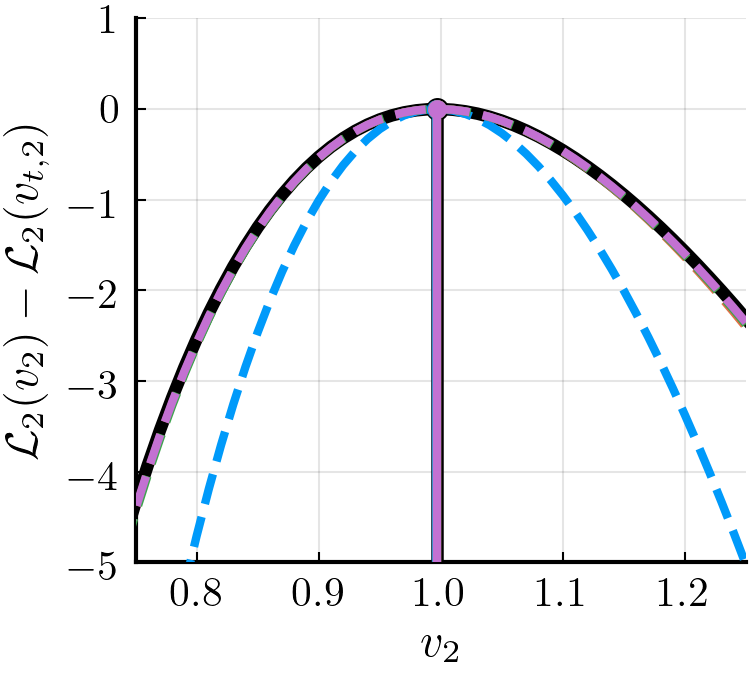}}
  
  \caption{%
    Convergence of alternating maximization w.r.t.~$\bmF$ and $\bmv$
    for various $\bmv$ updates.
    We consider $n = 10^3$ samples in $d = 10^2$ dimensions
    with $k=3$ underlying factors $\btllambda = (4,2,1)$.
    The noise is homoscedastic:
    both the first $n_1 = 200$ and remaining $n_2 = 800$ samples
    have noise variance $\tlv_1 = \tlv_2 = 1$.
    Walltimes are medians taken over 100 runs of the algorithm
    to reduce the effect of experimental noise.
    Markers denote each iteration.
  }
  \label{fig:comp:hom}
\end{figure*}

\begin{figure*} \centering
  \subfloat[Convergence w.r.t. $\clL$. \label{fig:comp:het:likelihood}]
    {\includegraphics[width=0.24\linewidth]{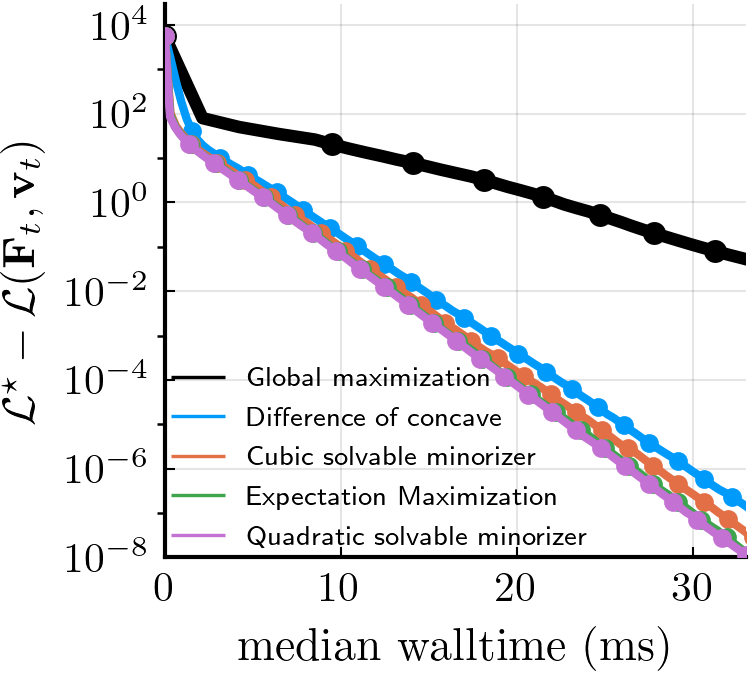}}%
  \hfill
  \subfloat[Convergence w.r.t. $\bmF$. \label{fig:comp:het:factors}]
    {\includegraphics[width=0.24\linewidth]{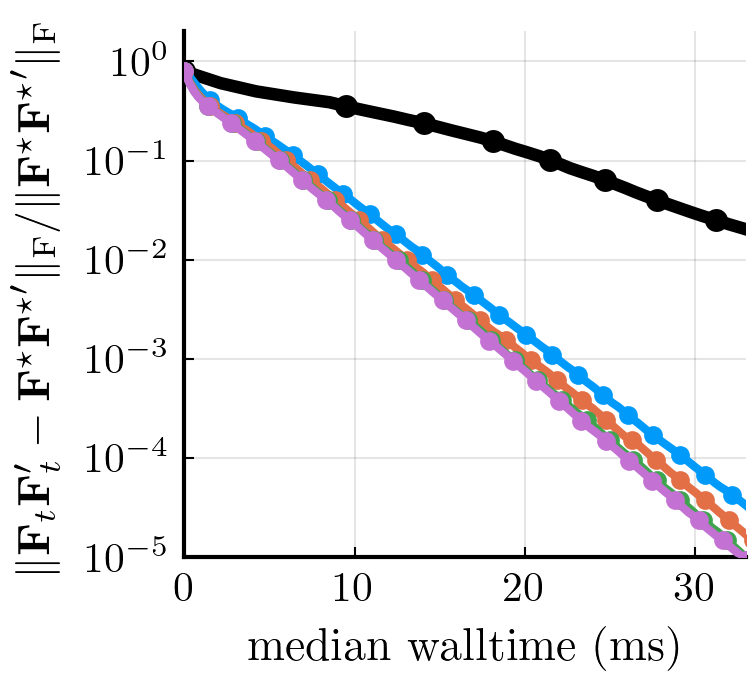}}%
  \hfill
  \subfloat[Minorizers w.r.t. $\bmv$ at homoscedastic PPCA initialization. \label{fig:comp:het:minorizer}]
    {\includegraphics[width=0.24\linewidth]{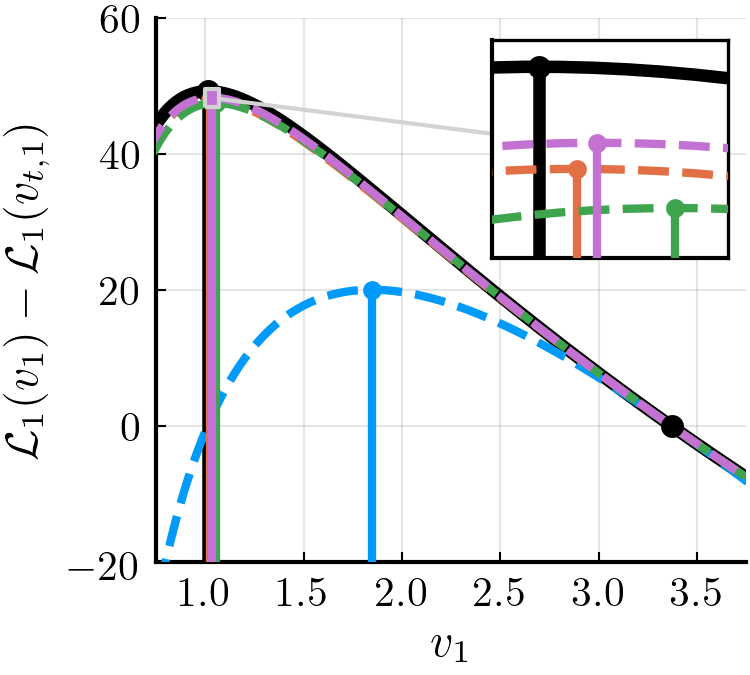}\;\;
     \includegraphics[width=0.24\linewidth]{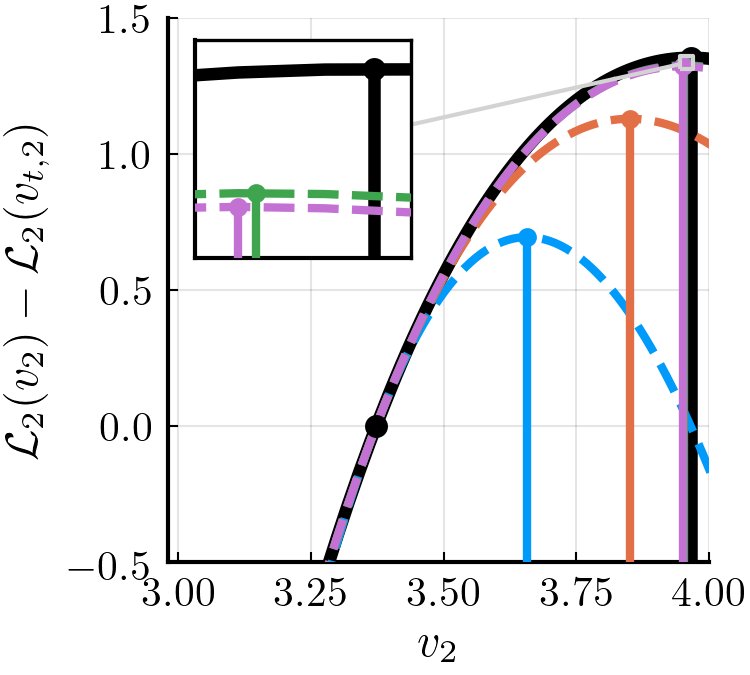}}
  
  \caption{%
    Same as \cref{fig:comp:hom} except here the noise is heteroscedastic
    with $\tlv_2 = 4$ and markers are placed every five iterations.
  }
  \label{fig:comp:het}
\end{figure*}

\Cref{sec:alg} described several update methods for \bmv
based on a variety of minorizers.
It is not obvious a priori which choice is best,
so this section compares their relative performance.
We consider $n = 10^3$ samples in $d = 10^2$ dimensions
generated according to the model \cref{eq:ppca:generative_model}
with $k=3$ factors generated
as $\btlF = \btlU\diag^{1/2}(\btllambda)$,
where $\btlU = (\btlu_1,\dots,\btlu_k) \in \bbR^{d \times k}$
is drawn uniformly at random from among
$d \times k$ matrices having orthonormal columns%
\footnote{Specifically, drawn according to the Haar measure on the Stiefel manifold,
  see, e.g., \cite[Section 2.5.1]{chikuse2003sos},
  as implemented in the Julia package
  \href{https://github.com/JuliaManifolds/Manifolds.jl}{Manifolds.jl}.}
and $\btllambda = (4,2,1)$.
The first $n_1 = 200$ samples have noise variance $\tlv_1 = 1$,
and the remaining $n_2 = 800$ have noise variance $\tlv_2$.

\Cref{fig:comp:hom} considers
the homoscedastic setting where $\tlv_2 = 1$
(yet the two variances are still unknown to the algorithm).
As a baseline, we take $\bmF^\star$ and $\bmv^\star$
to be the solutions obtained by 1000 iterations
of Expectation Maximization updates for both $\bmF$ and $\bmv$.
The associated log-likelihood is $\clL^\star \defequ \clL(\bmF^\star,\bmv^\star)$.
\Cref{fig:comp:hom:likelihood}
plots convergence of the log-likelihood $\clL^\star - \clL(\bmF_t,\bmv_t)$
versus walltime for iterates $(\bmF_t,\bmv_t)$
obtained by the various choices for the $\bmv$ update.
Note that $\clL^\star - \clL(\bmF_t,\bmv_t)$
is the log of the likelihood-ratio
between the converged solution $(\bmF^\star,\bmv^\star)$
and iterate $(\bmF_t,\bmv_t)$.
\Cref{fig:comp:hom:factors}
plots convergence for the $\bmF$ iterates
with respect to the normalized factor difference
$
\normfrobs{ \bmF_t \bmF_t' - \bmF^\star {\bmF^\star}' } / 
\normfrobs{ \bmF^\star {\bmF^\star}'}
$.
Iterations are indicated on both plots by the markers,
and walltime only includes the updates themselves
(i.e., not calculation of the log-likelihood).

Among the $\bmv$ update methods,
\rootfind
typically ascends the log-likelihood the most per iteration,
but is also the most computationally expensive.
As a result, it converges more slowly with respect to walltime.
The difference of concave update is computationally cheaper
but ascends the least per iteration initially.
The final three updates
(Cubic solvable MM, Expectation Maximization, Quadratic solvable MM)
have fairly similar computational cost
and log-likelihood increase per iteration.

The \rootfind update
corresponds to maximizing the univariate functions $\clL_\ell(v_\ell)$.
The remaining update methods each correspond to maximizing an associated minorizer.
\Cref{fig:comp:hom:minorizer} plots these minorizers
at the homoscedastic PPCA initialization (\cref{sec:init}),
shifted to be zero at the current iterate.
For this homoscedastic case,
the homoscedastic PPCA initialization is already close to optimal
and the minorizers (with the exception of the difference of concave minorizer)
closely follow the log-likelihood.

\Cref{fig:comp:het} considers a heteroscedastic case with $\tlv_2 = 4$.
As in the homoscedastic case,
\rootfind
converges the most slowly overall
due to its high computational cost per iteration
(more so in fact).
Likewise, the difference of concave update is again computationally cheaper
but ascends the least per iteration initially,
and the remaining three update methods converge the most rapidly.
\Cref{fig:comp:het:minorizer} illustrates the comparative tightness
of the various minorizers at the homoscedastic PPCA initialization.
The initialization is far from optimal for this heteroscedastic case,
and the relative differences in tightness among the minorizers are
more clearly visible.
Based on these experiments,
we recommend using the EM minorizer
or the quadratic solvable minorizer
for the \bmv updates.

\section{Statistical performance experiments}
\label{sec:exp:stat}

This section evaluates the statistical performance
of \heppca through simulation.
We consider $n = 10^3$ samples in $d = 10^2$ dimensions
generated according to the model \cref{eq:ppca:generative_model}
with $k=3$ factors generated
as $\btlF = \btlU\diag^{1/2}(\btllambda)$,
where $\btlU = [\btlu_1,\dots,\btlu_k] \in \bbR^{d \times k}$
is drawn uniformly at random from among
$d \times k$ matrices having orthonormal columns
and $\btllambda = (4,2,1)$.
The first $n_1 = 200$ samples have noise variance $\tlv_1 = 1$,
and the remaining $n_2 = 800$ have $\tlv_2 = \sigma_2^2$,
where we sweep $\sigma_2$ from $0$ to $3$.
We use $100$ iterations
of alternating EM updates for $\bmF$ and $\bmv$
with the homoscedastic PPCA initialization.

\begin{figure*} \centering
    \begin{tikzpicture}
    \legendentry{ 0.0}{0}{}      {homoppca}  {\footnotesize Homoscedastic PPCA: Full data}
    \legendentry{ 5.5}{0}{}      {group1}    {\footnotesize PPCA: Group 1}
    \legendentry{ 8.9}{0}{}      {group2}    {\footnotesize PPCA: Group 2}
    \legendentry{12.3}{0}{dashed}{heteroppca}{\footnotesize \heppca}
    \end{tikzpicture}\\[-4mm]
    \subfloat[Normalized factor estimation error \label{fig:stat:perf:homo:F}]
      {\includegraphics[width=0.24\linewidth]{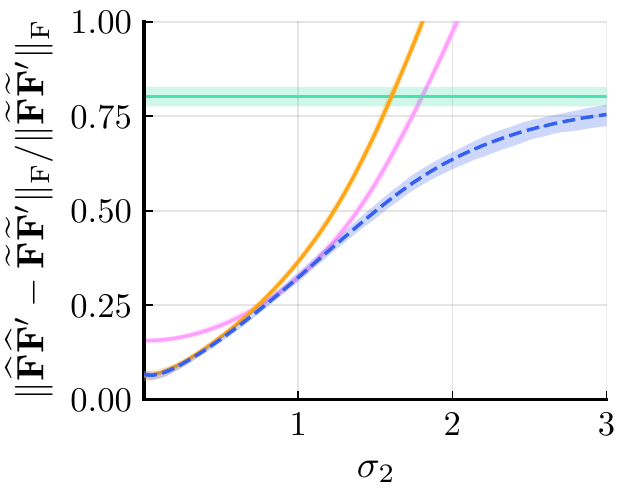}}
    \hfill
    \subfloat[Component 1 recovery \label{fig:stat:perf:homo:u1}]
      {\includegraphics[width=0.24\linewidth]{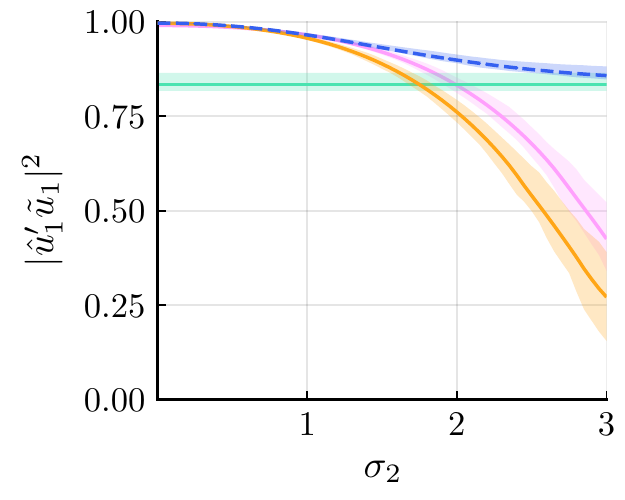}}
    \hfill
    \subfloat[Component 2 recovery \label{fig:stat:perf:homo:u2}]
      {\includegraphics[width=0.24\linewidth]{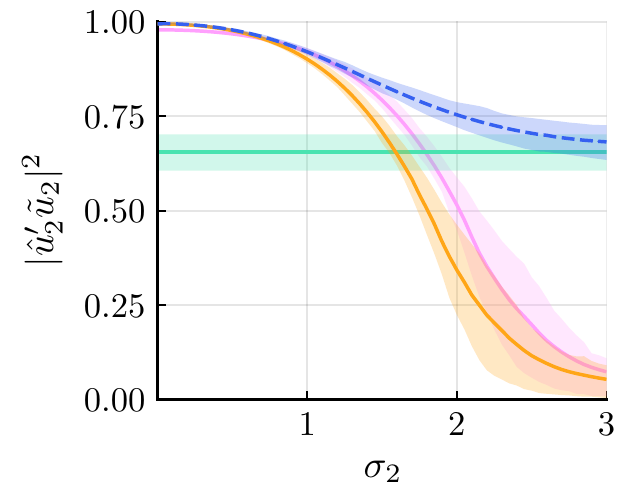}}
    \hfill
    \subfloat[Component 3 recovery \label{fig:stat:perf:homo:u3}]
      {\includegraphics[width=0.24\linewidth]{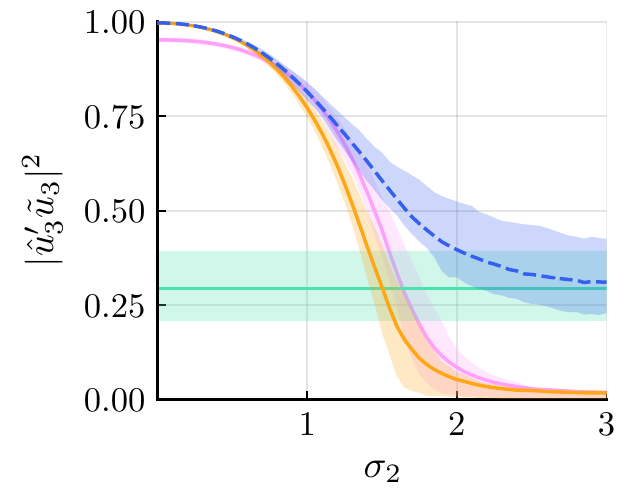}}
    \caption{%
      Comparison with homoscedastic PPCA applied on:
      i) full data,
      ii) only group 1,
      i.e, the $n_1 = 200$ samples with noise variance $\tlv_1=1$,
      and
      iii) only group 2,
      i.e., the $n_2=800$ samples with noise variance $\tlv_2=\sigma_2^2$.
      Lower is better in (a),
      and higher is better in (b)-(d).
      Heteroscedastic PPCA (\heppca) outperforms
      the homoscedastic methods
      on all four metrics.
      The mean and interquartile intervals (25th to 75th percentile) from $100$ data realizations
      are shown as curves and ribbons, respectively.
    }
    \label{fig:stat:perf:homo}
\end{figure*}

\subsection{Comparison with homoscedastic methods}
\label{sec:stat:perf:homo}

\Cref{fig:stat:perf:homo} compares
the recovery of the latent factors $\btlF$
by \heppca
with those obtained by applying homoscedastic PPCA on:
a) the full data,
b) only group 1,
i.e, the $n_1 = 200$ samples with noise variance $\tlv_1=1$,
and
c) only group 2,
i.e., the $n_2=800$ samples with noise variance $\tlv_2=\sigma_2^2$.
These homoscedastic PPCA approaches are reasonable and common choices
in the absence of reliable heteroscedastic algorithms;
it is worthwhile to understand their performance.
The mean and interquartile intervals (25th to 75th percentile) from $100$ data realizations
are shown as curves and ribbons, respectively.

\Cref{fig:stat:perf:homo:F} plots the normalized factor covariance estimation error,
defined as
$\normfrobs{\bhtF\bhtF' - \btlF\btlF'} / \normfrobs{\btlF\btlF'}$
where $\bhtF \in \bbR^{d \times k}$ is the estimated factor matrix.
\Cref{fig:stat:perf:homo:u1,fig:stat:perf:homo:u2,fig:stat:perf:homo:u3}
plot the component recoveries $|\bhtu_1'\btlu_1|^2,\dots,|\bhtu_3'\btlu_3|^2$,
where $\bhtu_1,\dots,\bhtu_3 \in \bbR^d$
are the principal eigenvectors of $\bhtF\bhtF'$.
Lower is better for estimation error and higher is better for component recovery.

When $\sigma_2$ is small enough,
homoscedastic PPCA applied
to only group 2 performs the best among
the homoscedastic PPCA's.
In this case, group 2 is relatively clean,
and the components are reliably recovered.
Using the full data incorporates more samples,
but in this case including the noisier group 1 data
does more harm than good
since homoscedastic PPCA treats them uniformly.
There is a tradeoff here
between having more samples
and including noisier samples.
Finally,
using only group 1 performs worst;
it is smaller and noisier.

With increasing $\sigma_2$,
the performance of homoscedastic PPCA degrades
when applied to the full data or group 2
since they incorporate these increasingly noisy samples.
The effect is more pronounced for using only group 2,
and eventually the tradeoff reverses;
using the full data becomes best among the homoscedastic PPCA options.
In particular,
when $\sigma_2 = 1$,
the full data actually has homoscedastic noise
and there is no statistical benefit
to using only either group.
As $\sigma_2$ continues to increase,
group 2 data is eventually so noisy
that it becomes best to only use group 1.
Just past $\sigma_2 > 1$, however,
using only group~1 remains worse than using only the noisier group~2 data.
In this regime,
the more abundant samples in group 2
win out
over the cleaner samples in group 1.
Which homoscedastic PPCA option performs best
depends crucially on the interplay
of these tradeoffs,
making it unclear a priori
which to use.

\heppca uses all the data
but estimates and accounts for the heteroscedastic noise.
In \cref{fig:stat:perf:homo},
it essentially matches or outperforms
all three homoscedastic PPCA options
across the entire range of $\sigma_2$.
In particular,
for small $\sigma_2$,
it closely matches the performance
of using only group 2,
and for large $\sigma_2$
it closely matches that of using only group 1.
In some sense,
it appears to automatically ignore unreliable data.
For moderate $\sigma_2$,
it outperforms the three homoscedastic PPCA options.
In this regime,
it is suboptimal to ignore either group of data
or to use both but treat them uniformly,
and \heppca appropriately combines them.

Notably, \heppca performs nearly the same across this sweep
as the variant developed in \cite{hong2019ppf}
that assumed \emph{known} noise variances,
even though the noise variances are now unknown and jointly estimated.
See \cref{fig:stat:perf:heppca:known}.
\Cref{sec:stat:perf:v_lambda,sec:stat:perf:blocks} study
the quality of the noise variance estimates.

\begin{figure*} \centering
  \begin{tikzpicture}
  \legendentry{ 0.0}{0}{}      {inv}            {\footnotesize Weighted PCA (inv. noise var.)}
  \legendentry{ 5.0}{0}{}      {sqinv}          {\footnotesize Weighted PCA (square inv. noise var.)}
  \legendentry{10.8}{0}{}      {heteropca}      {\footnotesize HeteroPCA \cite{zhang2018hpa:arxiv:v2}}
  \legendentry{14.0}{0}{dashed}{heteroppca}     {\footnotesize \heppca}
  \end{tikzpicture}\\[-4mm]
  \subfloat[Normalized subspace est. error \label{fig:stat:perf:hetero:U}]
    {\includegraphics[width=0.24\linewidth]{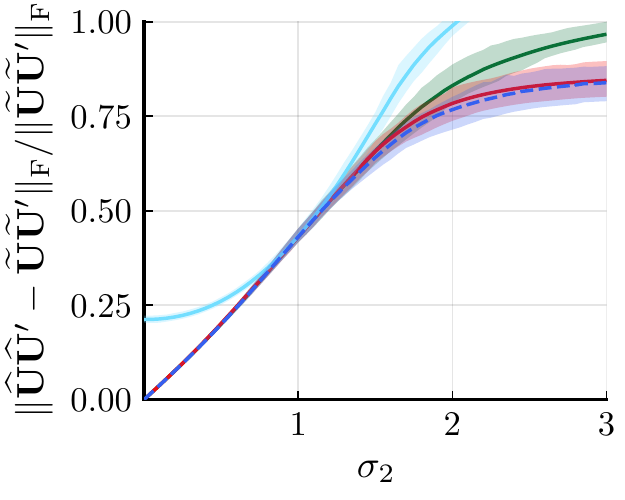}}
  \hfill
  \subfloat[Component 1 recovery \label{fig:stat:perf:hetero:u1}]
    {\includegraphics[width=0.24\linewidth]{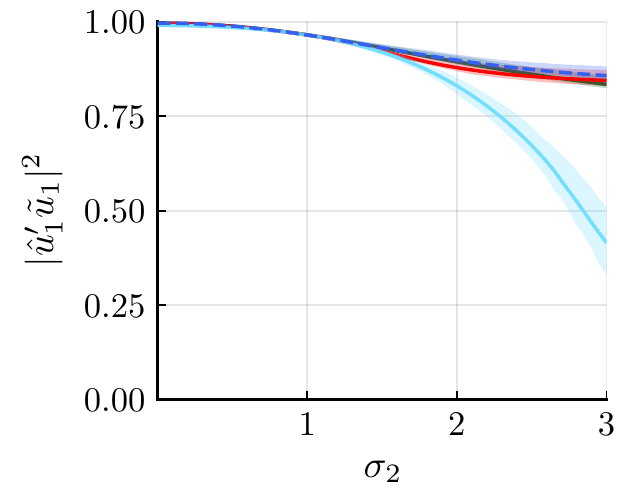}}
  \hfill
  \subfloat[Component 2 recovery \label{fig:stat:perf:hetero:u2}]
    {\includegraphics[width=0.24\linewidth]{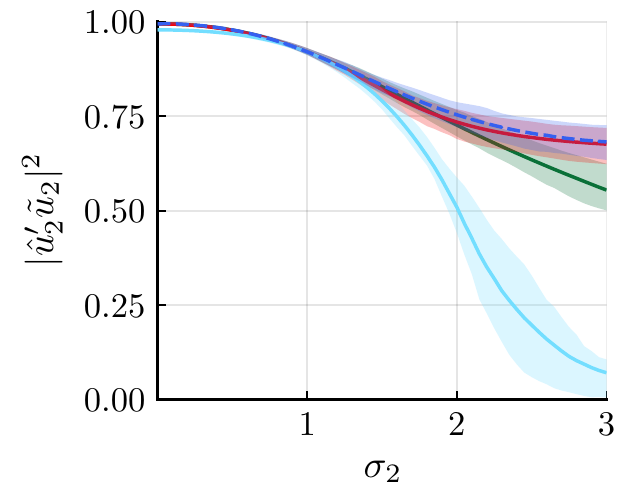}}
  \hfill
  \subfloat[Component 3 recovery \label{fig:stat:perf:hetero:u3}]
    {\includegraphics[width=0.24\linewidth]{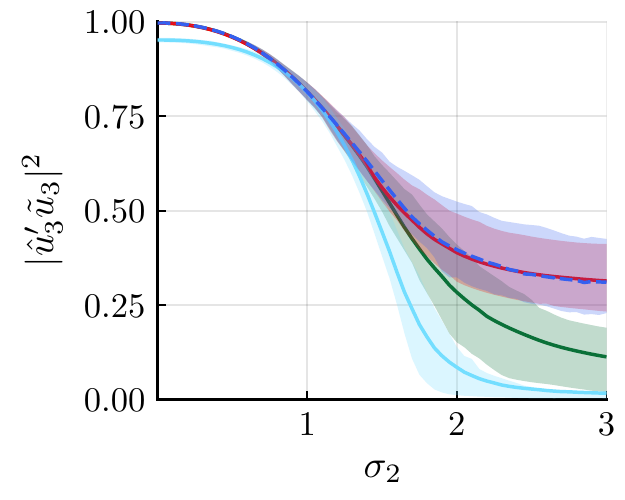}}
  \caption{%
    Comparison with heteroscedastic methods:
    HeteroPCA \cite{zhang2018hpa:arxiv:v2}
    and weighted PCA (inverse and square inverse noise variance weights
    calculated using the true noise variances).
    Lower is better in (a),
    and higher is better in (b)-(d).
    \heppca is among the best heteroscedastic methods.
    The mean and interquartile intervals (25th to 75th percentile) from $100$ data realizations
    are shown as curves and ribbons, respectively.
  }
  \label{fig:stat:perf:hetero}
\end{figure*}

\subsection{Comparison with heteroscedastic methods}
\label{sec:stat:perf:hetero}

\Cref{fig:stat:perf:hetero} compares
the recovery of the latent components $\btlU$
by \heppca
with those obtained by
HeteroPCA~\cite{zhang2018hpa:arxiv:v2}
and by weighted PCA
with:
a) inverse noise variance weights,
and b) square inverse noise variance weights.
HeteroPCA is an iterative method
designed for noise that is heteroscedastic
within each sample;
here we use $10$ iterations.
Weighted PCA is a simple and efficient variant of PCA
that accounts for heteroscedasticity across samples
by down-weighting noisier samples.
A typical choice for weights
is inverse noise variance weighting
that effectively rescales samples
so their (scaled) noise becomes homoscedastic.
Square inverse noise variance weighting
is a choice of weights
that more aggressively down-weights noisier samples.
It can be more effective for weak signals,
as revealed by the analysis in \cite{hong2018owp:arxiv:v2}.

\Cref{fig:stat:perf:hetero:U} plots the normalized subspace estimation error,
defined as
$\normfrobs{\bhtU\bhtU' - \btlU\btlU'} / \normfrobs{\btlU\btlU'}$
where we denote the estimated factor eigenvectors
as $\bhtU = [\bhtu_1,\dots,\bhtu_k] \in \bbR^{d \times k}$.
\Cref{fig:stat:perf:hetero:u1,fig:stat:perf:hetero:u2,fig:stat:perf:hetero:u3}
plot the corresponding component recoveries $|\bhtu_1'\btlu_1|^2,\dots,|\bhtu_3'\btlu_3|^2$.
As before,
the mean and interquartile intervals (25th to 75th percentile) from $100$ data realizations
are shown as curves and ribbons, respectively,
and lower is better for estimation error
and higher is better for component recovery.

When $\sigma_2$ is small,
both weighted PCA methods perform similarly to \heppca.
They both down-weight group 1 samples
and benefit from the clean group 2 samples.
As $\sigma_2$ increases,
a gap in the performance appears between
inverse noise variance and square inverse noise variance weighted PCA.
In this regime,
the more aggressive square inverse noise variance weights perform better
for the weaker second and third components.
For the stronger first component,
inverse noise variance weights remain comparable
and are, in fact, better for moderate $\sigma_2$.
Throughout the sweep,
\heppca matches
or slightly outperforms
the statistical performance
of both weighted PCA methods.
The relatively favorable performance
of these methods highlights
the benefit of accounting for heteroscedasticity across samples.
Doing so enables them to make better use of all the available data.
Note that unlike the two weighted PCA methods,
\heppca is not given the noise variances
and instead estimates them.
Moreover,
to apply the weighted PCA methods shown here,
one must choose between the two weights
(neither is uniformly better than the other),
whereas \heppca works well across the range of noise variances.

HeteroPCA also accounts for heteroscedastic noise,
but does so primarily for heteroscedasticity within each sample,
rather than across samples.
Heteroscedasticity within each sample
biases the diagonal of the covariance matrix,
even in expectation,
and HeteroPCA corrects for this.
However, it treats the samples themselves fairly uniformly.
Consequently,
its performance here closely resembles
that of homoscedastic PPCA on the full data,
as shown in \cref{fig:stat:perf:homo}.
This behavior highlights
the qualitative difference between
heteroscedasticity within and heteroscedasticity across samples;
they manifest differently and must both be addressed.
\Cref{sec:stat:perf:hetero:within_sample}
considers a setting with noise that is heteroscedastic in both ways
and illustrates the opportunity for further works
considering the combination.

\begin{figure} \centering
  \oneortwocol{
    \hfill
    \subfloat[Noise variance recovery \label{fig:stat:perf:v_lambda:v}]
      {\includegraphics[width=0.275\linewidth]{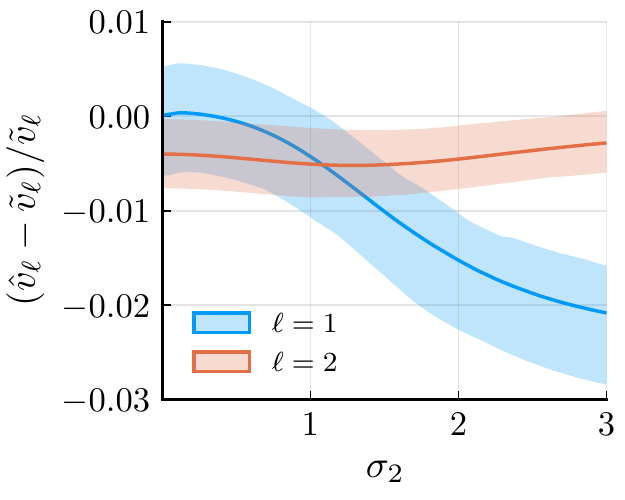}}
    \hfill
    \subfloat[Factor eigenvalue recovery \label{fig:stat:perf:v_lambda:lambda}]
      {\includegraphics[width=0.275\linewidth]{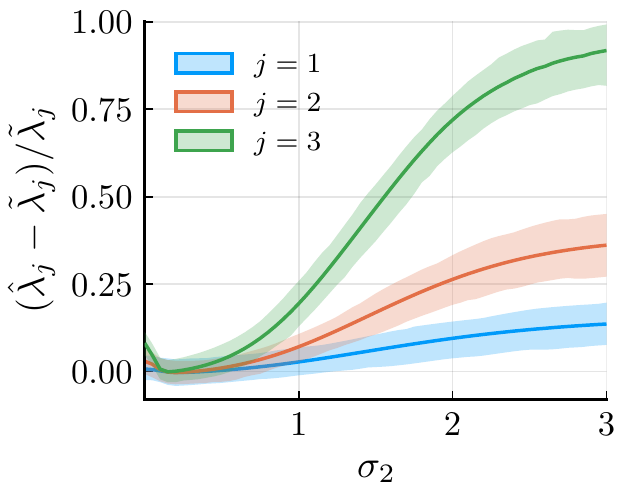}}
    \hfill
  }{
    \subfloat[Noise variance recovery for each group $\ell = 1,2$. \label{fig:stat:perf:v_lambda:v}]
      {\includegraphics[width=0.49\linewidth]{statperf,sim1,v}}
    \hfill
    \subfloat[Factor eigenvalue recovery for each component $j = 1,2,3$. \label{fig:stat:perf:v_lambda:lambda}]
      {\includegraphics[width=0.49\linewidth]{statperf,sim1,lambda}}
  }
  \caption{Relative bias of \heppca estimates for
    noise variances and factor eigenvalues.
    The mean and interquartile intervals (25th to 75th percentile) from $100$ data realizations
    are shown as curves and ribbons, respectively.
  }
  \label{fig:stat:perf:v_lambda}
\end{figure}

\subsection{Bias in estimated noise variances and factor eigenvalues}
\label{sec:stat:perf:v_lambda}

\Cref{fig:stat:perf:v_lambda} plots the relative biases
of the estimated noise variances $\htv_1$ and $\htv_2$,
as well as those of the estimated factor covariance eigenvalues
$\bhtlambda = (\htlambda_1,\dots,\htlambda_3)$.
As before,
the mean and interquartile intervals (25th to 75th percentile) from $100$ data realizations
are shown as curves and ribbons, respectively,
but now closer to zero is better.
Positive values mean that \heppca has overestimated,
and negative values indicate that it has underestimated.
Taken together, \cref{fig:stat:perf:v_lambda:v,fig:stat:perf:v_lambda:lambda}
show a general negative bias in the estimated noise variances
paired with a general positive bias in the estimated factor eigenvalues.
This behavior is consistent with
a corresponding behavior for homoscedastic PPCA
in the setting of homoscedastic noise \cite{passemier2015oeo}.
Providing a similar characterization
for \heppca
and a corresponding de-biasing procedure
is an exciting, but nontrivial,
direction for future work.

\subsection{Dependence of noise variance estimates on block sizes}
\label{sec:stat:perf:blocks}

\begin{figure} \centering
  \oneortwocol{
    \subfloat[%
      Noise variance estimates for a
      single data realization. \label{fig:stat:perf:blocks:v:single}]
      {\includegraphics[width=0.5\linewidth]{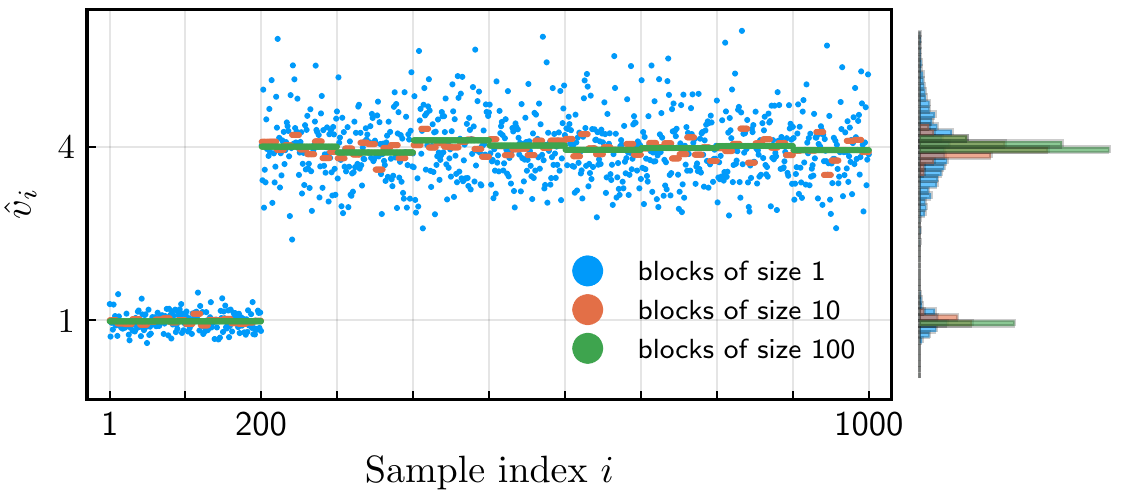}}
    \subfloat[%
      Mean and interquartile intervals from
      $100$ data realizations. \label{fig:stat:perf:blocks:v:100}]
      {\includegraphics[width=0.5\linewidth]{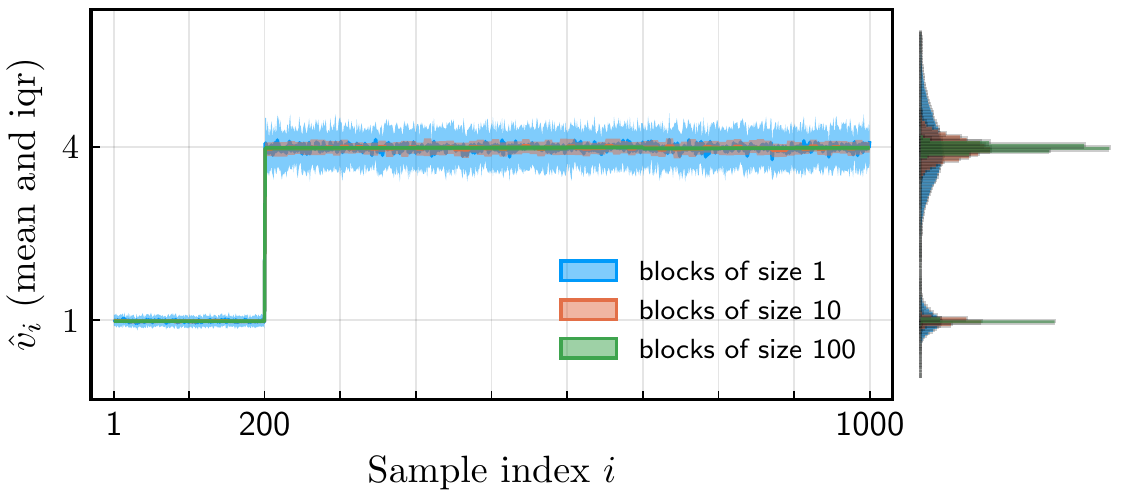}}
  }{
    \subfloat[%
      Noise variance estimates for a
      single data realization. \label{fig:stat:perf:blocks:v:single}]
      {\includegraphics[width=0.9\linewidth]{statperf,sim2,singlerun}}
    
    \subfloat[%
      Mean and interquartile intervals from
      $100$ data realizations. \label{fig:stat:perf:blocks:v:100}]
      {\includegraphics[width=0.9\linewidth]{statperf,sim2,allruns,100}}
  }
  \caption{Estimated noise variances for varying block sizes.}
  \label{fig:stat:perf:blocks:v}
\end{figure}

\Cref{fig:stat:perf:blocks:v}
fixes the noise variance of group 2 at $\tlv_2 = 4$,
i.e., $\sigma_2 = 2$,
and shows the noise variances estimated
for all of the $n=10^3$ samples
when the samples are passed to \heppca
in (non-overlapping) blocks of size 1, 10 and 100.
Doing so reveals how the estimates depend on block size,
and captures settings where the true latent groupings are unknown.
Notably, a block size of 1 incorporates no a priori knowledge of the groupings.
\Cref{fig:stat:perf:blocks:v:single}
shows a single representative data realization,
and
\cref{fig:stat:perf:blocks:v:100}
shows the mean and interquartile intervals
(25th to 75th percentile)
obtained from $100$ data realizations.
Both include corresponding histograms
on the right
showing the distributions of estimated noise variances.

\begin{figure} \centering
  \oneortwocol{
    \includegraphics[width=0.5\linewidth]{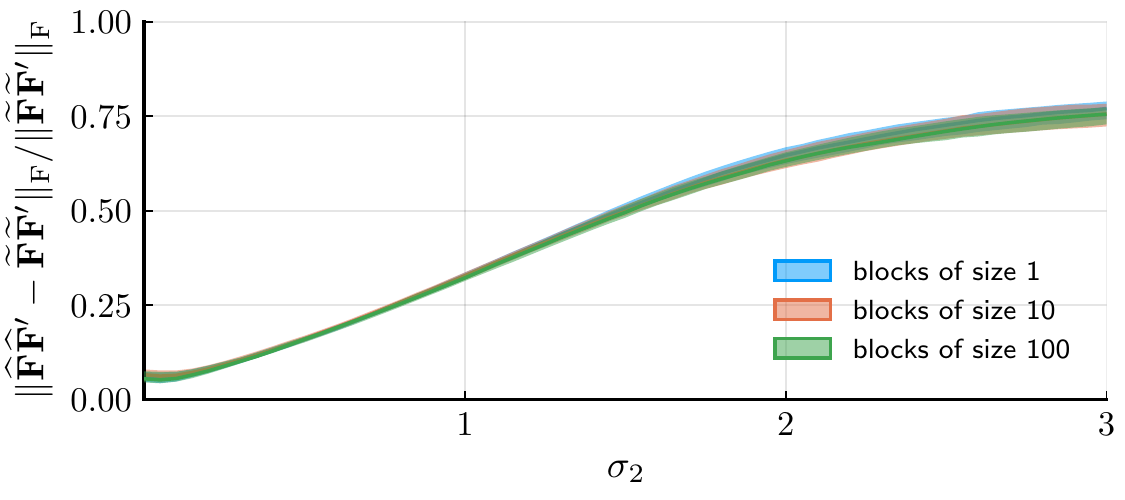}
  }{
    \includegraphics[width=0.9\linewidth]{statperf,sim3,F,med.pdf}
  }
  \caption{Normalized factor estimation error (median and interquartile intervals)
    for varying block sizes.
    The three block sizes have practically identical performance.}
  \label{fig:stat:perf:blocks:F}
\end{figure}

Notably, the estimates are fairly concentrated
around the true latent noise variances of $\btlv = (1,4)$
at blocks of size $100$
even though this choice splits
the first group of $n_1 = 200$ samples into two groups
and the second group of $n_2 = 800$ into eight groups.
These groups are visible in \cref{fig:stat:perf:blocks:v:single}
as bars that tie together samples in the same block.
Interestingly, blocks of size $10$ are not much more noisy,
while being significantly less restrictive.
Moreover, even using blocks of size $1$,
at which point each sample is allowed its own noise variance estimate,
provides relatively reliable estimates
that cluster around the latent noise variances.
The data contain enough information
to obtain reasonable estimates of these noise variances.
Nevertheless,
when samples can be reasonably grouped together into blocks,
e.g., grouping them by source or sensor,
doing so can significantly denoise the estimates
even when the blocks are relatively small.
An interesting direction for future work
is to jointly estimate these clusters
from the data.

\Cref{fig:stat:perf:blocks:F}
plots the corresponding normalized factor estimation errors
for $v_2 = \sigma_2^2$ where $\sigma_2$ ranges from $0$ to $3$.
The median and interquartile intervals (25th to 75th percentile) from $100$ data realizations
are shown as curves and ribbons, respectively.
We use the median here because
the means for runs of block size 1,
which are likely the most challenging,
appeared to be skewed by outliers.
Notably,
all three block sizes perform quite closely
to \heppca with known blocks
and outperform homoscedastic variants
(cf. \cref{fig:stat:perf:homo:F}).

\begin{figure} \centering
  \oneortwocol{
    \subfloat[Single data realization (heatmap). \label{fig:stat:perf:hetero:within_sample:heatmap}]
      {\includegraphics[width=0.5\linewidth]{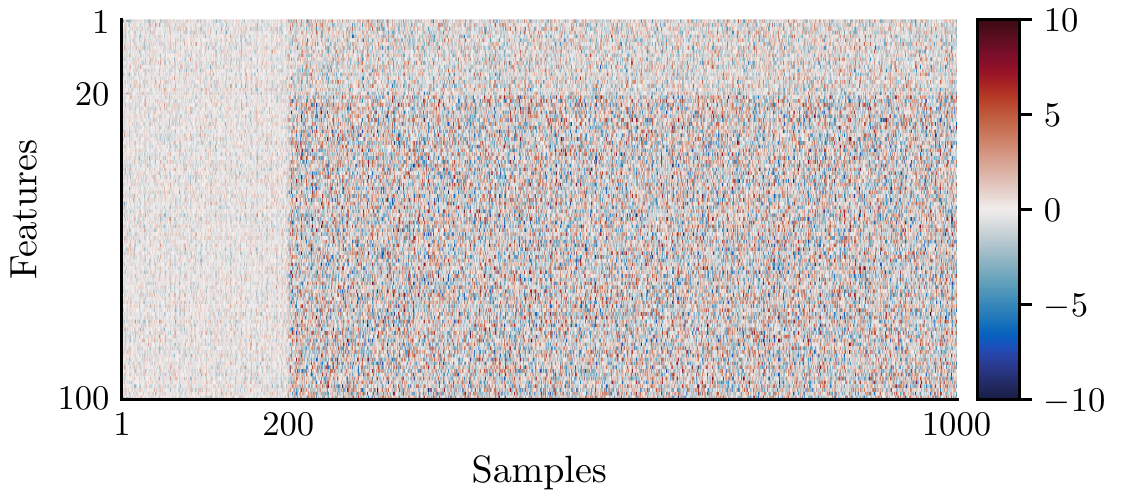}}
    \subfloat[Normalized subspace estimation error (lower is better). \label{fig:stat:perf:hetero:within_sample:U}]
      {\includegraphics[width=0.5\linewidth]{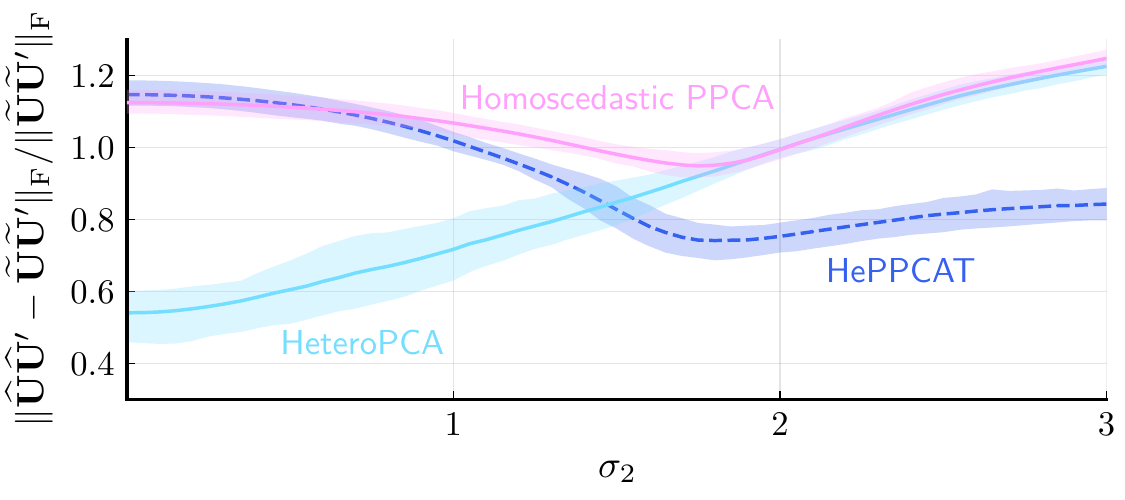}}
  }{
    \subfloat[Single data realization (heatmap). \label{fig:stat:perf:hetero:within_sample:heatmap}]
      {\includegraphics[width=0.9\linewidth]{statperf,sim4,data,heatmap}}
    
    \subfloat[Normalized subspace estimation error (lower is better). \label{fig:stat:perf:hetero:within_sample:U}]
      {\includegraphics[width=0.9\linewidth]{statperf,sim4,U}}
  }
  \caption{Heteroscedasticity across both features and samples.}
  \label{fig:stat:perf:hetero:within_sample}
\end{figure}

\subsection{Additional heteroscedasticity within samples}
\label{sec:stat:perf:hetero:within_sample}

\Cref{fig:stat:perf:hetero:within_sample}
considers data that is heteroscedastic
not just across samples,
but also within samples.
As before,
the first group of $n_1 = 200$ samples have noise variance $v_1 = 1$,
but now the second group of $n_2 = 800$ samples have
a noise variance fixed at $v_2^{(1)} = 4$
for the first $d^{(1)} = 20$ features
and noise variance $v_2^{(2)} = \sigma_2^2$
for the remaining $d^{(2)} = 80$ features,
where $\sigma_2$ ranges from $0$ to $3$.
Noise in the first group is homoscedastic within each sample,
but except for $\sigma_2 = 2$,
noise in the second group is heteroscedastic within each sample.
\Cref{fig:stat:perf:hetero:within_sample:heatmap}
shows an example data realization for $\sigma_2 = 3$;
observe that the first group are uniformly noisy,
the first 20 features of the second group are noisier,
and the final 80 features are noisiest.

\Cref{fig:stat:perf:hetero:within_sample:U}
plots the subspace estimation error across this range of heteroscedastic settings
for homoscedastic PPCA;
\heppca,
which accounts for heteroscedasticity across samples;
and HeteroPCA~\cite{zhang2018hpa:arxiv:v2},
which primarily
accounts for heteroscedasticity within each sample.
Namely, HeteroPCA accounts for bias in the diagonal of the covariance
that is caused by within-sample heteroscedasticity,
but treats the samples themselves uniformly.
When $\sigma_2$ is small,
accounting for heteroscedasticity within each sample
is more important and HeteroPCA is better.
However, the tradeoff reverses
and \heppca becomes better
as $\sigma_2$ grows towards $\sigma_2 = 2$,
at which point samples are heteroscedastic only across samples.
Interestingly,
heteroscedasticity across samples appears to continue to dominate
for $\sigma_2 > 2$.
Homoscedastic PPCA is generally worst,
as it does not account for either heteroscedasticity.
\heppca performs similarly for small $\sigma_2$,
where within-sample heteroscedasticity seems to dominate,
and HeteroPCA is similar for large $\sigma_2$,
where across-sample heteroscedasticity seems to dominate.
These results highlight a qualitative difference
between across-sample and within-sample heteroscedasticity;
both must be addressed.
Developing methods that simultaneously handle both types of heteroscedasticity,
which outperform all three methods across this range,
is an exciting direction for future work.

\section{Real Data Experiments}
\label{sec:exp:real_data}

\newcommand{\Ytest}{\bmY^{\mathrm{(test)}}}
\newcommand{\Ytrain}{\bmY^{\mathrm{(train)}}}
\newcommand{\ntest}{n^{\mathrm{(test)}}}
\newcommand{\ntrain}{n^{\mathrm{(train)}}}

This section applies HePPCAT to environmental monitoring data
containing
air quality measurements from both a few high precision instruments
and a large network of low-cost consumer-grade sensors.
High precision measurements are provided by
the U.S. Environmental Protection Agency (EPA) and its partners.
They maintain a nationwide network of Air Quality Index (AQI) sensor stations
that measure, monitor, and distribute air quality data
on the AirNow platform \cite{airnowEPA}.
The recent proliferation of low-cost consumer-grade AQI sensors,
such as PurpleAir \cite{purpleair}, provides a second source of data.
These sensors stream data continuously to developer platforms, such as ThingSpeak \cite{thingspeak},
creating a network of crowd-sourced air quality data
with greater spatial coverage and resolution
but generally lower precision.

\begin{figure}
    \centering
    \includegraphics[width=0.35\textwidth]{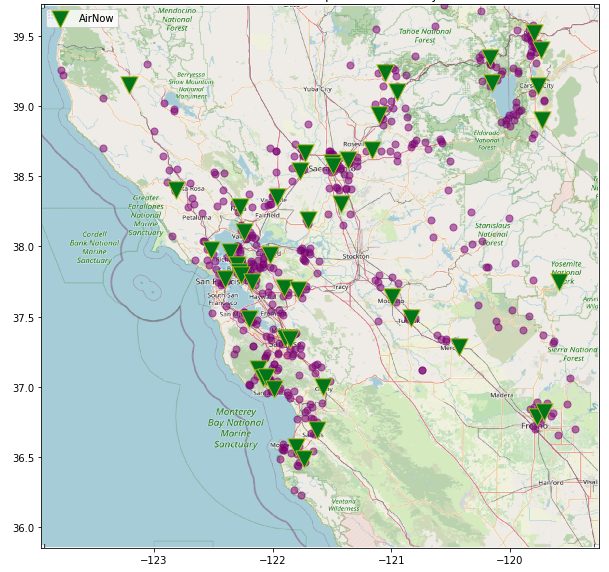}
    \caption{Sensor locations of AirNow (green triangles) and PurpleAir (purple circles) in the central California region of our experiments.}
    \label{fig:central_CA_map}
\end{figure}

We consider $\text{PM}_{2.5}$ particulate concentration readings (in $\mu g/m^3$)
from the AirNow platform and from outdoor PurpleAir sensors
across the central California region,
i.e., within longitudes (-123.948, -119.246) and latitudes (35.853, 39.724),
at the top of every hour from February 9-13, 2021.
We chose 10 random PurpleAir sensors nearby each of 46 AirNow sensors
to obtain balanced sensing coverage,
and omitted hours where at least one of the AirNow sensors did not record a measurement.
This gave $n_1=46$ AirNow samples $\bmy_{1,i}$ and $n_2=460$ PurpleAir samples $\bmy_{2,i}$,
where each sample $\bmy_{\ell,i}$ is a vector of $d=108$ readings of $\text{PM}_{2.5}$ across time.
\Cref{fig:central_CA_map} displays the map of the sensor locations for visualization.
AirNow measurements are calibrated and averaged over hour-long windows by the U.S. EPA,
whereas we collected the instantaneous readings from the PurpleAir sensors nearest to each hour.
We centered the two sensor groups $\bmY_1$ and $\bmY_2$ separately by subtracting from each its sample mean.

Lacking ground truth,
we evaluate how well the subspace learned by HePPCAT
on a subset of the samples generalizes to the rest.
Namely, we randomly select
$\ntrain_1 = 25$ AirNow samples and $\ntrain_2 = 250$ PurpleAir samples
as training data $\Ytrain \in \bbR^{d \times \ntrain}$.
The remaining
$\ntest_1 = 21$ AirNow samples and $\ntest_2 = 210$ PurpleAir samples
serve as test data $\Ytest \in \bbR^{d \times \ntest}$.
The training data $\Ytrain$ is then used
to estimate a basis for a $k=30$ dimensional subspace $\bhtU \in \bbR^{d \times k}$
using PPCA and HePPCAT.
We also consider PPCA-AN (PPCA on only the AirNow group $\Ytrain_1$)
and PPCA-PA (PPCA on only the PurpleAir group $\Ytrain_2$).
For each estimated $\bhtU$,
the performance on test data $\Ytest$
is quantified by the normalized root mean-squared error (NRMSE)
of the subspace reconstruction,
i.e.,
$\normfrobs{ \Ytest - \bhtU \bhtU' \Ytest } / \normfrob{\Ytest}$.
We also consider the corresponding NRMSE evaluated
on only the AirNow test data $\Ytest_1$
and on only the PurpleAir test data $\Ytest_2$,
as well as all the training data counterparts.

\begin{figure}
    \centering
    \hfill
    \subfloat[A single representative trial. \label{fig:aqi_heppcat_variances:rep}]
    {\includegraphics[width=1.61in]{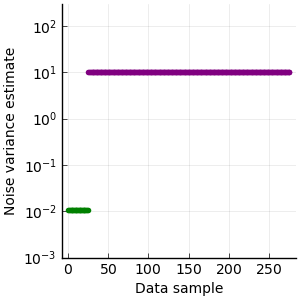}}
    \hfill
    \subfloat[Boxplots from all 200 trials. \label{fig:aqi_heppcat_variances:boxplot}]
    {\parbox[b]{1.6in}{\includegraphics[width=1.5in]{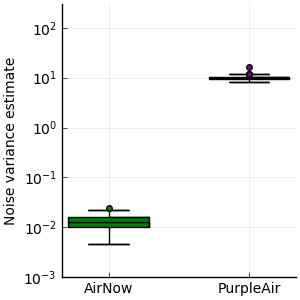}\vspace{2.8mm}}}
    \hfill
    \strut
    \caption{Noise variance estimates from HePPCAT
        for a single representative trial \protect\subref{fig:aqi_heppcat_variances:rep},
        and across all 200 trials \protect\subref{fig:aqi_heppcat_variances:boxplot}.
        In \protect\subref{fig:aqi_heppcat_variances:rep},
        the first 25 samples in green are from AirNow,
        and the remaining purple samples are from PurpleAir.
        A single noise variance is estimated for each group.
        Boxplots in \protect\subref{fig:aqi_heppcat_variances:boxplot}
        show the spread of these estimates across the 200 trials. Units are in $(\mu g / m^3)^2$.}
    \label{fig:aqi_heppcat_variances}
\end{figure}

We repeated this experiment for 200 random train-test splits of the data.
\Cref{fig:aqi_heppcat_variances:rep}
shows the noise variance estimates from HePPCAT
from a representative trial.
The estimated noise variance for the PurpleAir samples
is substantially higher than that for AirNow samples,
illustrating heterogeneity within this data.
This is reasonable given that the PurpleAir data
comes from low-cost consumer-grade sensors,
while the AirNow data comes from high precision instruments.
The corresponding box plots (indicating median, interquartile range, and outliers)
in \Cref{fig:aqi_heppcat_variances:boxplot}
show the spread of these estimates across the 200 trials.
The estimates remain fairly consistent.

\begin{figure*}
    \centering
    \subfloat[NRMSE with respect to full training data $\Ytrain$,
        AirNow training data $\Ytrain_1$,
        and PurpleAir training data $\Ytrain_2$.
        \label{fig:aqi_nrmse_box_plots:train}]
    {\includegraphics[width=0.99\linewidth]
    {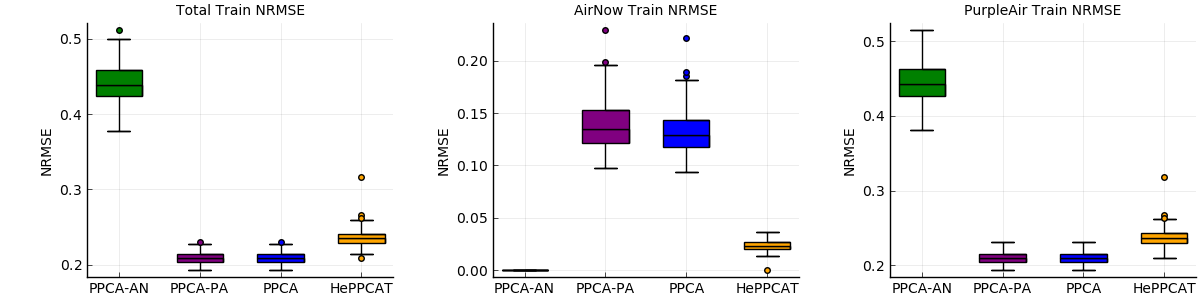}}
    
    \subfloat[NRMSE with respect to full test data $\Ytest$,
        AirNow test data $\Ytest_1$,
        and PurpleAir test data $\Ytest_2$.
        \label{fig:aqi_nrmse_box_plots:test}]
    {\includegraphics[width=0.99\linewidth]
    {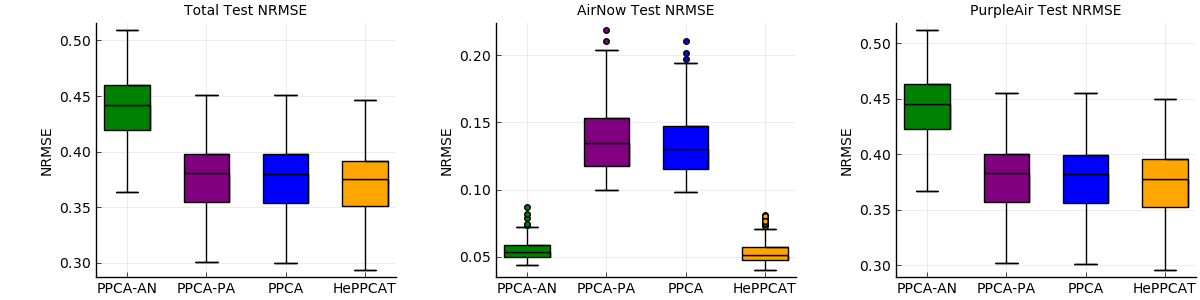}}
   
    \caption{Air quality data were split into training and test sets. A $k=30$ dimensional subspace basis $\bhtU$ was estimated from the training data using PPCA and HePPCAT, as well as PPCA-AN (PPCA trained using only the AirNow group)
    and PPCA-PA (PPCA trained using only the PurpleAir group). We evaluate the NRMSE $\normfrobs{ \bmY - \bhtU \bhtU' \bmY } / \normfrob{\bmY}$
    with respect to both training and test data (as well as their AirNow and PurpleAir subsets).
    Lower is better for all plots. As expected, HePPCAT is never best on training data. However, it is among the best on all test cases, indicating that it has found explanatory components across both data sources.}
    \label{fig:aqi_nrmse_box_plots}
\end{figure*}

\Cref{fig:aqi_nrmse_box_plots} shows boxplots for the training and test NRMSEs across the 200 trials.
As expected, PPCA (on full data) generally has the lowest training NRMSE on full data,
PPCA-AN is generally best on AirNow training data,
and likewise for PPCA-PA on PurpleAir training data.
Compared with PPCA (on full data),
HePPCAT has worse training NRMSE on the full data and on the PurpleAir data,
but has better training NRMSE on the cleaner AirNow data.
Turning to test NRMSE, however,
HePPCAT is among the best with respect to
not only AirNow data,
but \emph{also} the full data and even the PurpleAir data.
HePPCAT appears to more effectively leverage information
from both the cleaner but fewer AirNow samples
and the noisier but more numerous PurpleAir samples.

Overall,
the experiments here with air quality data
illustrate heterogeneity arising naturally in real data
and the potential for improved generalization by using HePPCAT.

\begin{figure*} \centering
    \oneortwocol{
        \includegraphics[width=0.8\textwidth]{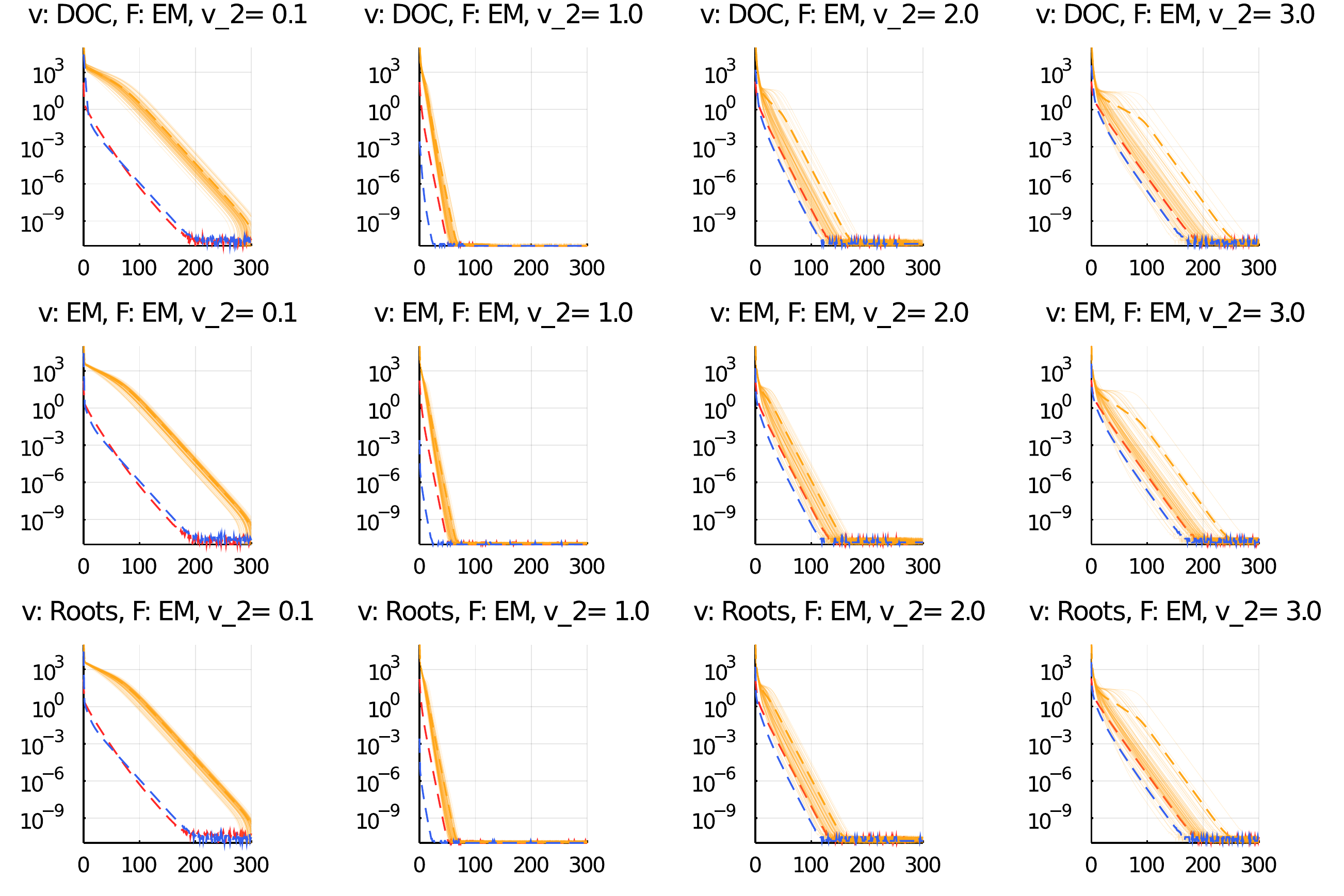} \label{fig:convergence_gaps_trace}
    }{
        \includegraphics[width=0.8\textwidth]{Fv_algs-convergence_gaps} \label{fig:convergence_gaps_trace}
    }
    
    \oneortwocol{
    \begin{tikzpicture}
    \begin{customlegend}[legend columns=-1]
    \addlegendimage{red,dashed}
    \addlegendentry{Oracle init (planted model)}
    \addlegendimage{blue,dashed}
    \addlegendentry{Homoscedastic PPCA init}
    \addlegendimage{orange,dashed}
    \addlegendentry{Random inits}
    \end{customlegend}
    \end{tikzpicture}
    }{
    \begin{tikzpicture}
    \begin{customlegend}[legend columns=-1]
    \addlegendimage{red,dashed}
    \addlegendentry{Oracle initialization (planted model)}
    \addlegendimage{blue,dashed}
    \addlegendentry{Homoscedastic PPCA initialization}
    \addlegendimage{orange,dashed}
    \addlegendentry{Random initializations}
    \end{customlegend}
    \end{tikzpicture}
    }

    \caption{Convergence gaps of each algorithm to the maximum converged log-likelihood per heteroscedastic noise experiment. $n = [200,800]$ and $\bmv_1 = 1.0$.}
    \label{fig:landscape_convergence_gaps}
\end{figure*}

\section{Investigation of the landscape}
\label{sec:exp:init}

This section empirically illustrates
the favorable landscape of our optimization problem
and our algorithms' convergence to the globally optimal solution
for synthetically generated data.
Even though the nonconvexity of the problem might lead one to wonder
if the choice of initialization matters,
we find that does not appear to be the case.
We generate the same low-rank heteroscedastic model as described in \cref{sec:exp:stat}
and sweep across $\sigma_2^2$ values 0.1, 1.0, 2.0, and 3.0.
The first regime should see \heppca largely down-weight group 1
and perform PCA on just group 2.
At $\sigma_2^2=1$, the dataset is statistically homoscedastic,
and we expect the landscape to behave similarly to that of PCA,
which enjoys a well-known landscape that features
no spurious local maxima and strict saddles \cite{ge2017no}.
As $\sigma_2^2$ increases, the distribution of the noise variances becomes more bimodal,
and the PPCA solution deviates farther from the optimal log-likelihood value.
In the low-noise end, the second data block has four times as many samples as the first block and a tenth of the noise variance.
In the noisiest setting,
the second data block has three times the noise variance.

For each noise setting, we record each algorithms' log-likelihood at iteration $t$ and in \cref{fig:landscape_convergence_gaps} show the difference
to the maximum log-likelihood found among all algorithms and trials.
We run 100 trials of each \heppca algorithm
with the initial estimate of $\bmF$ drawn randomly
with i.i.d. Gaussian $\clN(0,1)$ entries
and the initial estimate of $\bmv$ drawn randomly
with i.i.d. entries uniform on $[0,1)$.
We also examine initializations from the homoscedastic PPCA solution
and from the oracle planted model parameters.
The converged likelihood values for each algorithm and choice of initialization
concentrate tightly around the same maximum,
with this behavior consistent across a wide range of heteroscedastic noise levels,
indicating a well-behaved landscape.

In the homoscedastic regime,
the results are consistent with our expectation
that the PPCA initialization should be close to optimal,
as shown in the second column of \cref{fig:landscape_convergence_gaps}.
We observe as the data noise variances become more imbalanced,
the PPCA initialization becomes farther away from the global maximum,
but is still orders of magnitude better in likelihood
than random initialization.
The oracle initialization has the best likelihood
for all the heteroscedastic settings as expected,
but is still suboptimal since we are maximizing
a finite-sample likelihood.
An interesting direction is to study
how heteroscedastic noise
affects the likelihood and its maximum
in finite sample settings.


\section{Conclusion}
\label{sec:conc}

This paper developed efficient algorithms
for jointly estimating latent factors and noise variances
from data with heteroscedastic noise.
Maximizing the likelihood is a nontrivial nonconvex optimization problem,
and unlike the homoscedastic setting,
it is seemingly unsolvable via singular value decomposition.
The proposed algorithms alternate between updating the factor estimates
and the noise variance estimates,
with several choices for the noise variance update.
It is unclear a priori which choice is best,
and we compared their empirical convergence speeds in practice.
Further numerical experiments studying the statistical performance
highlighted the significant benefits
of properly accounting for heteroscedasticity.
Experiments on air quality data illustrated
heterogeneity arising naturally in real data
and improved generalization by using HePPCAT.
Given the nonconvexity of the problem,
one might wonder if initializing differently
could lead to better maximizers.
We provided empirical evidence that this is not the case;
the landscape, while nonconvex, appears favorable.

Extensions of the approach to handle more general settings,
e.g., missing data or additional heterogeneity across features,
are interesting directions for further work.
Likewise, there are many variations of PCA,
e.g., nonnegative matrix factorization,
and generalizations,
e.g., unions of subspaces,
that one could consider.
An extension to consider kernel PCA
would be interesting
\cite{tipping:01:skp,lawrence:05:pnl},
as noted by a reviewer.
One might also incorporate a clustering step
in the alternating algorithm
to estimate not only the noise variances
but also the blocks sharing a common noise variance.
Alternatively, one could consider the $L$ groupings
to be another latent variable in the log-likelihood,
and attempt to jointly estimate them.
Estimating the rank is another direction for further work.
Many classical methods were designed for homoscedastic noise,
and recent works,
e.g., \cite{leeb2018oss:arxiv:v4,hong2020stn:arxiv:v1,ke2020eot:arxiv:v1,landa2021brt:arxiv:v1},
have begun to explore this problem under heteroscedastic settings.
Some avenues for improving convergence speed
are
using momentum / extrapolation
\cite{chun:20:cao} of the alternating maximization updates,
as well as incremental variants.
One could also consider tackling the M-step in \cref{sec:em}
via an inner block coordinate ascent
with updates similar to those in \cref{sec:alg:F:em,sec:alg:v:em}
that ascend the EM minorizer \cref{eq:em:estep}.
This paper also raises several natural conjectures
about the landscape of the nonconvex objective,
which are beyond our present scope
and are exciting areas for further theoretical analysis.
Finally, it was observed in the homoscedastic setting
that noise variance estimates
tend to have a downward bias
that can be characterized and accounted for \cite{passemier2015oeo}.
A similar bias in variance estimates
appears to occur in the heteroscedastic setting,
and extending the previous approaches is a promising direction.

\appendix

\section{Challenges in Expectation Maximization}
\label{sec:em:challenges}

To carry out Expectation Maximization
via \cref{sec:em},
one might attempt to maximize \cref{eq:em:estep}
with respect to both $\bmv$ and $\bmF$
by first completing the square:
\oneortwocol{%
\begin{equation*}
\cbrL(\bmF, \bmv; \bmF_t, \bmv_t)
=
- \sum_{\ell=1}^L \bigg(
    \frac{d n_\ell}{2}\ln v_\ell + \frac{\normfrob{\bmY_\ell}^2}{2v_\ell}
  \bigg)
+ \frac{1}{2} \tr\Big\{
    \bmT(\bmv) \bmS(\bmv)^{-1} \bmT(\bmv)'
  \Big\}
-\frac{1}{2}
  \normfrob{
    \bmF \bmS(\bmv)^{1/2} - \bmT(\bmv) \bmS(\bmv)^{-1/2}
  }^2
,
\end{equation*}%
}{%
\begin{align*}
&\cbrL(\bmF, \bmv; \bmF_t, \bmv_t)
=
- \sum_{\ell=1}^L \bigg(
    \frac{d n_\ell}{2}\ln v_\ell + \frac{\normfrob{\bmY_\ell}^2}{2v_\ell}
  \bigg)
\\&\quad\nonumber
+ \frac{1}{2} \tr\Big\{
    \bmT(\bmv) \bmS(\bmv)^{-1} \bmT(\bmv)'
  \Big\}
\\&\quad\nonumber
-\frac{1}{2}
  \normfrob{
    \bmF \bmS(\bmv)^{1/2} - \bmT(\bmv) \bmS(\bmv)^{-1/2}
  }^2
,
\end{align*}%
}%
where the first two lines are constant with respect to $\bmF$, and
\begin{align*}
\bmT(\bmv) &\defequ \sum_{\ell=1}^L \frac{1}{v_\ell} \bmY_\ell \bbrZ_{t,\ell}'
, \oneortwocol{&}{\\}
\bmS(\bmv) &\defequ
  \sum_{\ell=1}^L \frac{1}{v_\ell}
    (\bbrZ_{t,\ell}\bbrZ_{t,\ell}' + n_\ell v_{t,\ell} \bmM_{t,\ell})
.
\end{align*}
Thus \cref{eq:em:estep} is maximized
with respect to $\bmF$ by
$\bmF = \bmT(\bmv)\bmS(\bmv)^{-1}$,
yielding the maximization problem with respect to $\bmv$ of
\begin{align} \label{eq:em:optF:sub}
&\cbrL(\bmT(\bmv)\bmS(\bmv)^{-1}, \bmv; \bmF_t, \bmv_t)
\oneortwocol{}{\\&\nonumber}
=
\frac{1}{2} \tr\Big\{
  \bmT(\bmv) \bmS(\bmv)^{-1} \bmT(\bmv)'
\Big\}
-
\sum_{\ell=1}^L \bigg(
  \frac{d n_\ell}{2}\ln v_\ell + \frac{\normfrob{\bmY_\ell}^2}{2v_\ell}
\bigg)
.
\end{align}
\Cref{eq:em:optF:sub} is not easily optimized with respect to $\bmv$
for $L > 1$
because of the matrix product in the trace.
This term introduces coupling among the noise variances $v_\ell$
that may keep the problem from separating into $L$ univariate problems.
Intuitively, the noise variances $v_\ell$ appear to be coupled via
their impact on the optimal latent factors $\bmF$.
Novel approaches to efficiently optimize \cref{eq:em:optF:sub},
e.g., by studying its critical points,
is an interesting direction for future work.


\section{Comparison with known noise variances}

\begin{figure} \centering
  \begin{tikzpicture}
    \legendentry{0.0}{0}{dashed}{heteroppca}     {\footnotesize \heppca}
    \legendentry{2.5}{0}{}      {heteroppcaknown}{\footnotesize Oracle variant using known \bmv \cite{hong2019ppf}}
  \end{tikzpicture}\\[-4mm]
  \oneortwocol{
    \subfloat[Normalized subspace est. error \label{fig:stat:perf:heppca:known:U}]
      {\includegraphics[width=0.24\linewidth]{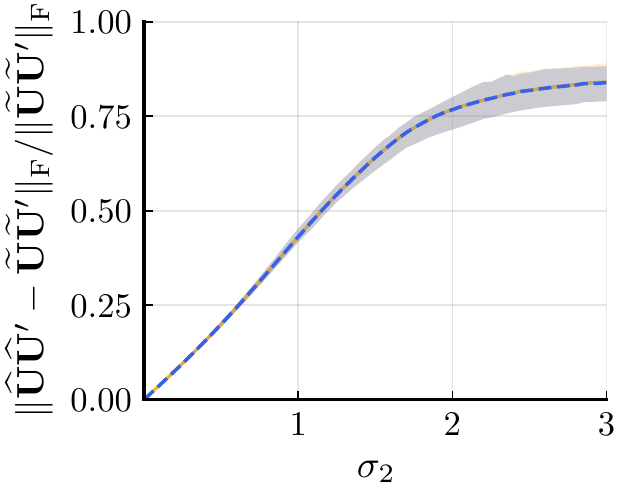}}
    \hfill
    \subfloat[Component 1 recovery \label{fig:stat:perf:heppca:known:u1}]
      {\includegraphics[width=0.24\linewidth]{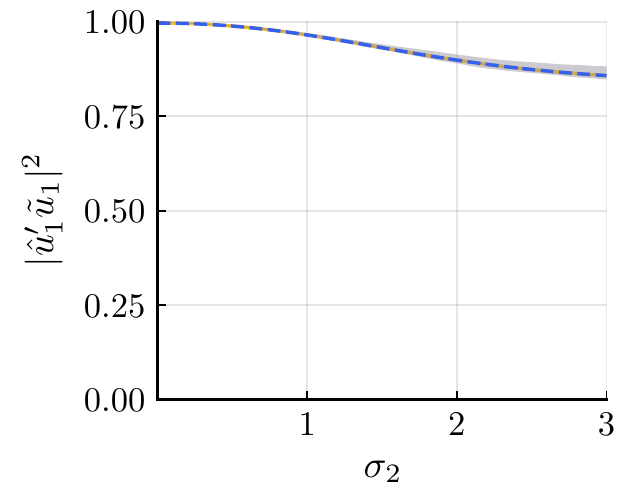}}
    \hfill
    \subfloat[Component 2 recovery \label{fig:stat:perf:heppca:known:u2}]
      {\includegraphics[width=0.24\linewidth]{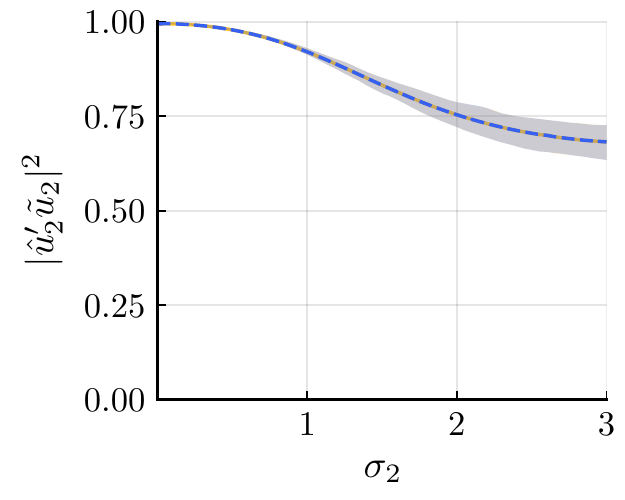}}
    \hfill
    \subfloat[Component 3 recovery \label{fig:stat:perf:heppca:known:u3}]
      {\includegraphics[width=0.24\linewidth]{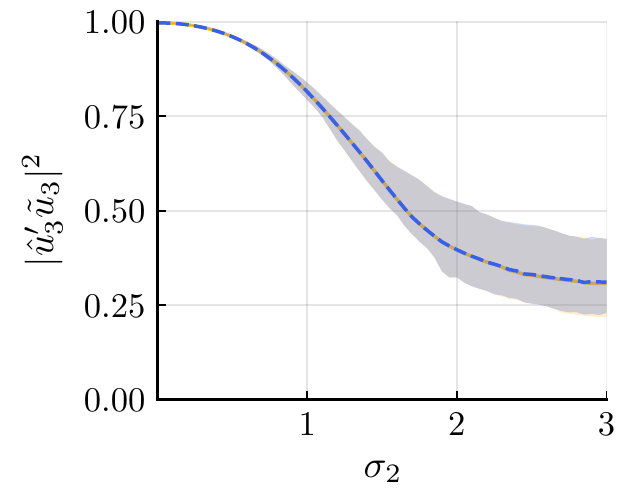}}
  }{
    \subfloat[Normalized subspace est. error \label{fig:stat:perf:heppca:known:U}]
      {\includegraphics[width=0.48\linewidth]{statperf,sim1c,U}}
    \hfill
    \subfloat[Component 1 recovery \label{fig:stat:perf:heppca:known:u1}]
      {\includegraphics[width=0.48\linewidth]{statperf,sim1c,u1}}
    \hfill
    \subfloat[Component 2 recovery \label{fig:stat:perf:heppca:known:u2}]
      {\includegraphics[width=0.48\linewidth]{statperf,sim1c,u2}}
    \hfill
    \subfloat[Component 3 recovery \label{fig:stat:perf:heppca:known:u3}]
      {\includegraphics[width=0.48\linewidth]{statperf,sim1c,u3}}
  }
  \caption{%
    Comparison of the statistical performance of \heppca
    with the variant developed in \cite{hong2019ppf}
    that assumed \emph{known} noise variances,
    under the experimental sweep in \cref{sec:stat:perf:homo,sec:stat:perf:hetero}.
    Lower is better in (a),
    and higher is better in (b)-(d).
    \heppca performs nearly the same
    even though the noise variances are now unknown and jointly estimated.
  }
  \label{fig:stat:perf:heppca:known}
\end{figure}

\Cref{fig:stat:perf:heppca:known} 
compares \heppca with
the oracle variant developed in \cite{hong2019ppf}
that used \emph{known} noise variances;
the experimental setup matches that of \cref{sec:stat:perf:homo,sec:stat:perf:hetero}.
Even though the noise variances are unknown and jointly estimated in \heppca,
the performance is nearly the same.

\section*{Acknowledgments}

The authors thank Arnaud Breloy for helpful discussions,
pointing us to applications in RADAR,
and sharing useful references to relevant works.
They thank Hoon Hong for helpful discussions
related to finding roots of high-degree univariate polynomials.
They thank Raj Nadakuditi for helpful suggestions on how to evaluate the performance of HePPCAT on the air quality data.
They thank the anonymous reviewers for their feedback and suggestions.

\bibliographystyle{IEEEtran}
\bibliography{refs}

\end{document}